\documentclass[moor]{informs1}              

\usepackage[english]{babel}
\usepackage{amsmath}
\usepackage{amstext}
\usepackage{amssymb}
\usepackage{mathtools}
\usepackage{natbib}
 \NatBibNumeric
 \def\bibfont{\small}%
 \bibpunct[, ]{[}{]}{,}{n}{}{,}%

\usepackage[colorlinks=true,breaklinks=true,bookmarks=true,urlcolor=blue,
     citecolor=blue,linkcolor=blue,bookmarksopen=false,draft=false]{hyperref}

\def\EMAIL#1{\href{mailto:#1}{#1}}

\TheoremsNumberedBySection  

\EquationsNumberedBySection 

\def\rF{\mathbb{F}}
\def\R{\mathbb{R}}

\def\diam{\mathop{\rm diam}}
\def\Lip{\mathop{\rm Lip}}

\def\Pr{\mathop{\rm Pr}}
\def\intr{\mathop{\rm int}}

\def\argmin{\mathop{\rm arg\, min}}

\def\cE{\mathbb{E}}

\def\hJ{\hat{J}}

\def\hf{\hat{f}}

\def\B{{\mathcal B}}

\def\P{{\mathcal P}}

\def\S{{\mathcal S}}

\def\sX{{\mathsf X}}

\def\sA{{\mathsf A}}
\def\sH{{\mathsf H}}
\def\sZ{{\mathsf Z}}
\def\sM{{\mathsf M}}
\def\sR{{\mathsf R}}
\def\sE{{\mathsf E}}

\begin{document}

\RUNAUTHOR{Saldi, Y\"uksel, and Linder}
\RUNTITLE{Asymptotic Optimality of Finite Approximations to MDPs}

\TITLE{On the Asymptotic Optimality of Finite Approximations to Markov Decision Processes with Borel Spaces}

\ARTICLEAUTHORS{
\AUTHOR{Naci Saldi}
\AFF{Coordinated Science Laboratory, University of Illinois,Urbana, IL 61801-2307, USA.\\
\{\EMAIL{nsaldi@illinois.edu}\}}
\AUTHOR{Serdar Y\"uksel, Tam\'{a}s Linder}
\AFF{Department of Mathematics and Statistics, Queen's University, Kingston, ON, Canada, K7L 3N6. \{\EMAIL{yuksel,linder@mast.queensu.ca}\}}}

\ABSTRACT{
Calculating optimal policies is known to be computationally difficult for Markov decision processes (MDPs) with Borel state and action spaces. This paper studies finite-state approximations of discrete time Markov decision processes with Borel state and action spaces, for both discounted and average costs criteria. The stationary policies thus obtained are shown to approximate the optimal stationary policy with arbitrary precision under quite general conditions for discounted cost and more restrictive conditions for average cost. For compact-state MDPs, we obtain explicit rate of convergence bounds quantifying how the approximation improves as the size of the approximating finite state space increases. Using information theoretic arguments, the order optimality of the obtained convergence rates is established for a large class of problems. We also show that, as a pre-processing step the action space can also be finitely approximated with sufficiently large number points; thereby, well known algorithms, such as value or policy iteration, Q-learning, etc., can be used to calculate near optimal policies.
}

\MSCCLASS{93E20, 90C40, 90C39}
\ORMSCLASS{Primary: Dynamic programming/optimal control, probability ; secondary: Infinite state, Markov processes}

\maketitle

\section{Introduction.}
\label{sec0}

In this paper, our goal is to study the finite-state approximation problem for computing near optimal policies for discrete time Markov decision processes (MDPs) with Borel state and action spaces, under discounted and average costs criteria. Although the existence and structural properties of optimal policies have been studied extensively in the literature, computing such policies is generally a challenging problem for systems with uncountable state spaces. This situation also arises in the fully observed reduction of a partially observed Markov decision process even when the original system has finite state and action spaces (see, e.g., \citet{YuBe04}).

As has been extensively studied in the literature (see, e.g., \citet{chow1991optimal} and the literature review below), one way to compute approximately optimal solutions for such MDPs is to construct a reduced model with a new transition probability and a one-stage cost function by quantizing the state/action spaces, i.e., by discretizing them on a finite grid. We exhibit that under quite general continuity conditions on the one-stage cost function and the transition probability for the discounted cost and under some additional restrictions on the ergodicity properties of Markov chains induced by deterministic stationary policies for the average cost, the optimal policy for the approximating finite model applied to the original model has cost that converges to the optimal cost, as the discretization becomes finer. Moreover, under additional continuity conditions on the transition probability and the one stage cost function we also obtain bounds for a rate of approximation in terms of the number of points used to discretize the state space, thereby providing a tradeoff between the computation cost and the performance loss in the system. In particular, we study the following two problems.
\begin{itemize}
\item[\textbf{(Q1)}] Under what conditions on the components of the MDP do the true costs corresponding to the optimal policies obtained from finite models converge to the optimal value function as the number of grid points goes to infinity? For this problem, we are only concerned with the convergence of the approximation; that is, we do not establish bounds for a rate of approximation.
\item[\textbf{(Q2)}] Can we obtain explicit bounds on the performance loss due to the discretization in terms of the number of grid points if we strengthen the  conditions sufficient~in~\textbf{(Q1)}?
\end{itemize}

Combined with our recent works \citet{SaLiYu13-2,SaYuLi16}, where we investigated the asymptotic optimality of the quantization of action sets, the results in this paper lead to a constructive algorithm for obtaining approximately optimal solutions. First the action space is quantized with small error, and then the state space is quantized with small error, which results in a finite model that well approximates the original MDP. When the state space is compact, we also obtain rates of convergence for both approximations, and using information theoretic tools we establish that the obtained rates of convergence are order-optimal for a given class of MDPs. Since there exist various computational algorithms for finite-state Markov decision problems, the analysis in this paper can be considered to be {\it constructive}.

Various methods have been developed to compute approximate value functions and near optimal policies. A partial list of these techniques is as follows: approximate dynamic programming, approximate value or policy iteration, simulation-based techniques, neuro-dynamic programming (or reinforcement learning), state aggregation, etc. For rather complete surveys of these techniques, we refer the reader to \citet{Fox71,Whi78,Whi79,Lan81,BeTs96,ReKr02,Ort07,Whi80,Whi82,Ber75,DuPr13,DuPr14} and references therein. With the exception of \citet{DuPr14,Ort07}, these papers in general study either the finite horizon cost or the discounted infinite horizon cost. Also, the majority of these results are for MDPs with discrete (i.e., finite or countable) state and action spaces, or a bounded one-stage cost function (e.g., \citet{Fox71,Whi78,Whi79,Roy06,Whi80,Whi82,Cav86,BeTs96,ReKr02,Ort07,Ber75}). Those that consider general state and action spaces (see, e.g., \citet{DuPr12,DuPr13,DuPr14,Ber75,chow1991optimal}) assume in general Lipschitz type continuity conditions on the components of the control model, in order to provide a rate of convergence analysis for the approximation error. Some of the results only consider approximating the value function and do not provide a procedure to compute near optimal policies (e.g., \citet{Lan81,Whi79,DuPr13}).

Our paper differs from these results in the following ways: (i) we consider a general setup,
where the state and action spaces are Borel (with the action space being compact), and the one-stage cost function is possibly unbounded, (ii) since we do not aim to provide rate of convergence result in the first problem \textbf{(Q1)}, the continuity assumptions we impose on the components of the control model are weaker than the conditions imposed in prior works that considered general state and action spaces, (iii) we also consider the challenging average
cost criterion under reasonable assumptions. The price we pay for imposing weaker assumptions in \textbf{(Q1)} is that we do not obtain explicit performance bounds in terms of the number of grid points used in the approximations. However, such bounds can be obtained under further assumptions on the transition probability and the one-stage cost functions; this is considered in problem \textbf{(Q2)} for compact-state MDPs.

Our approach to solve problem \textbf{(Q1)} can be summarized as follows: (i) first, we obtain approximation results for the compact-state case, (ii) we find conditions under which a compact representation leads to near optimality for non-compact state MDPs, (iii) we prove the convergence of the finite-state models to non-compact models. As a by-product of this analysis, we obtain {\it compact-state-space approximation}s for an MDP with non-compact Borel state space. In particular, our findings directly lead to finite models if the state space is countable; similar problems in the countable context have been studied in the literature for the discounted cost; see \citet[Section 6.10.2]{Put05}.

We note that the proposed method for solving the approximation problem for compact-state MDPs with the discounted cost is partly inspired by \citet{Roy06}. Specifically, we generalize the operator proposed for an approximate value iteration algorithm in \citet{Roy06} to uncountable state spaces. Then, unlike in \citet{Roy06}, we use this operator as a transition step between the original optimality operator and the optimality operator of the approximate model. In \citet{Ort07}, a similar construction was given for finite state-action MDPs. Our method to obtain finite-state MDPs from the compact-state model can be regarded as a generalization of this construction. We note that a related work of \citet{DuPr14} develops a sequence of approximations using empirical distributions of an underlying probability measure with respect to which the transition probability of the MDP is absolutely continuous. By imposing Lipschitz type continuity conditions on the components of the control model, \citet{DuPr14} obtains a concentration inequality type upper bound on the accuracy of the approximation based on the Wasserstein distance of order 1 between the probability measure and its empirical estimate. These conditions are stronger than what we impose for the problem \textbf{(Q1)}. We note that \citet{DuPr14} adopts a simulation based approximation leading to probabilistic guarantees on the approximation, whereas we adopt a quantization based approach leading to deterministic approximation guarantees. For a review of further simulation based methods, see e.g., \citet{ChFuHuMa07,JaVa06}.

The approach developed in the paper is also useful in networked control applications where transmission of real-valued actions to an actuator is not realistic when there is an information transmission constraint between a plant, a controller, and an actuator (see, e.g., \citet{YuBa13}). On the other hand, the elements of a finite action set can be transmitted across a finite capacity information channel. Even though the problem of optimal quantization for information transmission from a plant/sensor to a controller has been studied extensively (see, e.g. references in \citet{YuBa13}), these type of results appear to be new in the networked control literature when the problem of transmitting signals from a controller to an actuator is considered. Furthermore, tools from information theory allow for obtaining lower bounds on the approximation performance; using such an argument we show that the construction in this paper is order-optimal for a large class of models.

The rest of the paper is organized as follows. In Section~\ref{compact} we study the approximation problem \textbf{(Q1)} for MDPs with compact state space. In Section~\ref{sec2} an analogous approximation result is obtained for MDPs with non-compact state space. Discretization of the action space is considered in Section~\ref{act dist} for a general state space. In Section~\ref{compact:rateconv} we derive quantitative bounds on the approximation error in terms of the number of points used to discretize the state space for the compact-state case. In Section~\ref{compact:order} the order optimality of the obtained bounds on the approximation errors is established.
In Section~\ref{examples} we present an example to numerically illustrate our results. Section~\ref{conc} concludes the paper.

\subsection{Notation and Conventions.}\label{intro notation}

For a metric space $\sE$, the Borel $\sigma$-algebra (the smallest $\sigma$-algebra that contains the open sets of $\sE$) is denoted by $\B(\sE)$.  We let $B(\sE)$ and $C_b(\sE)$ denote the set of all bounded Borel measurable and continuous real functions on $\sE$, respectively. For any $u \in C_b(\sE)$ or $ u \in B(\sE)$, let $\|u\| \coloneqq \sup_{e \in \sE} |u(e)|$
which turns $C_{b}(\sE)$ and $B(\sE)$ into Banach spaces. Given any Borel measurable function $w: \sE \rightarrow [1,\infty)$ and any real valued Borel measurable function $u$ on $\sE$, we define the $w$-norm of $u$ as
\begin{align}
\|u\|_w \coloneqq \sup_{e\in\sE} \frac{|u(e)|}{w(e)}, \nonumber
\end{align}
and let $B_w(\sE)$ denote the Banach space of all real valued measurable functions $u$ on $\sE$ with finite $w$-norm; see \citet{HeLa99}. Let $\P(\sE)$ denote the set of all probability measures on $\sE$. A sequence $\{\mu_n\}$ of probability measures on $\sE$ is said to converge weakly (resp., setwise) (see \citet{HeLa03}) to a probability measure $\mu$ if $\int_{\sE} g(e) \mu_n(de)\rightarrow \int_{\sE} g(e) \mu(de) \text{ for all } g\in C_b(\sE)$ (resp., for all $g \in B(\sE)$). For any $\mu,\nu \in \P(\sE)$, the total variation distance between $\mu$ and $\nu$, denoted as $\|\mu - \nu\|_{TV}$, is equivalently defined as
\begin{align}
\|\mu-\nu\|_{TV} &\coloneqq 2 \sup_{D \in \B(\sE)} |\mu(D) - \nu(D)| = \sup_{\|g\| \leq 1} \biggl| \int_{\sE} g(e) \mu(de) - \int_{\sE} g(e) \nu(de) \biggr|. \nonumber
\end{align}
Unless otherwise specified, the term `measurable' will refer to Borel measurability in the rest of the paper.

\subsection{Markov Decision Processes.}\label{sec1}

A discrete-time Markov decision process (MDP) can be described by a five-tuple
\begin{align}
\bigl( \sX, \sA, \{\sA(x): x \in \sX\}, p, c \bigr), \nonumber
\end{align}
where Borel spaces (i.e., Borel subsets of complete and separable metric spaces) $\sX$ and $\sA$ denote the \emph{state} and \emph{action} spaces, respectively. The collection $\{\sA(x): x \in \sX\}$ is a family of nonempty subsets $\sA(x)$ of $\sA$, which give the admissible actions for the state $x\in\sX$. The \emph{stochastic kernel} $p(\,\cdot\,|x,a)$ denotes the \emph{transition probability} of the next state given that previous state-action pair is $(x,a)$; see \citet{HeLa96}. Hence, it satisfies: (i) $p(\,\cdot\,|x,a)$ is an element of $\P(\sX)$ for all $(x,a)$, and (ii) $p(D|\,\cdot\,,\,\cdot\,)$ is a measurable function from $\sX\times\sA$ to $[0,1]$ for each $D\in\B(\sX)$. The \emph{one-stage cost} function $c$ is a measurable function from $\sX \times \sA$ to $\R$. In this paper, it is assumed that $\sA(x) = \sA$ for all $x\in\sX$.

Define the history spaces $\sH_0 = \sX$ and
$\sH_{t}=(\sX\times\sA)^{t}\times\sX$, $t=1,2,\ldots$ endowed with their
product Borel $\sigma$-algebras generated by $\B(\sX)$ and $\B(\sA)$. A
\emph{policy} is a sequence $\pi=\{\pi_{t}\}$ of stochastic kernels
on $\sA$ given $\sH_{t}$. The set of all policies is denoted by $\Pi$.
Let $\Phi$ denote the set of stochastic kernels $\varphi$ on $\sA$ given $\sX$, and let $\rF$ denote the set of all measurable functions $f$ from $\sX$ to $\sA$. A \emph{randomized Markov} policy is a sequence $\pi=\{\pi_{t}\}$ of stochastic kernels on $\sA$ given $\sX$. A \emph{deterministic Markov} policy is a sequence of stochastic kernels $\pi=\{\pi_{t}\}$ on $\sA$ given $\sX$ such that $\pi_{t}(\,\cdot\,|x)=\delta_{f_t(x)}(\,\cdot\,)$ for some $f_t \in \mathbb{F}$, where $\delta_z$ denotes the point mass at $z$. The set of randomized and deterministic Markov policies are denoted by $\sR\sM$ and $\sM$, respectively. A \emph{randomized stationary} policy is a
constant sequence $\pi=\{\pi_{t}\}$ of stochastic kernels on $\sA$ given $\sX$ such that
$\pi_{t}(\,\cdot\,|x)=\varphi(\,\cdot\,|x)$ for all $t$ for some
$\varphi \in \Phi$. A \emph{deterministic stationary} policy is a constant sequence of stochastic kernels $\pi=\{\pi_{t}\}$ on $\sA$ given $\sX$ such that $\pi_{t}(\,\cdot\,|x)=\delta_{f(x)}(\,\cdot\,)$ for all $t$ for some
$f \in \mathbb{F}$. The set of randomized and deterministic stationary policies are identified with the sets $\Phi$ and $\mathbb{F}$, respectively.

According to the Ionescu Tulcea theorem (see \citet{HeLa96}), an initial distribution $\mu$ on $\sX$ and a policy $\pi$ define a unique probability measure $P_{\mu}^{\pi}$ on $\sH_{\infty}=(\sX\times\sA)^{\infty}$.
The expectation with respect to $P_{\mu}^{\pi}$ is denoted by $\cE_{\mu}^{\pi}$.
If $\mu=\delta_x$, we write $P_{x}^{\pi}$ and $\cE_{x}^{\pi}$ instead of $P_{\delta_x}^{\pi}$ and $\cE_{\delta_x}^{\pi}$. The cost functions to be minimized in this paper are the $\beta$-discounted cost and the average cost, respectively given by
\begin{align}
J(\pi,x) &= \cE_{x}^{\pi}\biggl[\sum_{t=0}^{\infty}\beta^{t}c(x_{t},a_{t})\biggr], \nonumber \\
V(\pi,x) &= \limsup_{T\rightarrow\infty} \frac{1}{T} \cE_{x}^{\pi} \biggl[\sum_{t=0}^{T-1} c(x_t,a_t) \biggr]. \nonumber
\end{align}
With this notation, the discounted and average value functions of the control problem are defined as
\begin{align}
J^*(x) &\coloneqq \inf_{\pi \in \Pi} J(\pi,x), \nonumber \\
V^*(x) &\coloneqq \inf_{\pi \in \Pi} V(\pi,x). \nonumber
\end{align}
A policy $\pi^{*}$ is said to be optimal if $J(\pi^{*},x) = J^*(x)$ (or $V(\pi^{*},x) = V^*(x)$ for the average cost) for all $x \in \sX$. Under fairly mild conditions, the set $\rF$ of deterministic stationary policies contains an optimal policy
for discounted cost (see, e.g., \citet{HeLa96,FeKaZa12}) and average cost optimal control problems (under somewhat stronger continuity/recurrence conditions, see, e.g., \citet{FeKaZa12}).

\begin{remark}
We note that the path-wise infinite sum $\sum_{t=0}^{\infty} \beta^t c(x_t, a_t)$ may not be well-defined in the definition of $J$ if $c$ is only assumed to be measurable. However, further assumptions that will be imposed in later sections ensure that $J$ is a well-defined function.
\end{remark}

\subsection{Auxiliary Results}

To avoid measurability problems associated with the operators that will be defined for the approximation problem in the discounted cost case, it is necessary to enlarge the set of functions on which these operators can act. To this end, in this section we review the notion of analytic sets and lower semi-analytic functions, and state the main results that will be used in the sequel to tackle these measurability problems. For a detailed treatment of analytic sets and lower semi-analytic functions, we refer the reader to \citet{ShBe79,BlFrOr74}, \citet[Chapter 39]{Kur66}, and \citet[Chapter 7]{BeSh78}.

Let $\mathbb{N}^{\infty}$ be the set of sequences of natural numbers endowed with the product topology. With this topology, $\mathbb{N}^{\infty}$ is a complete and separable metric space. A subset $A$ of a Borel space $\sE$ is said to be \emph{analytic} if it is a continuous image of $\mathbb{N}^{\infty}$. Note that Borel sets are always analytic.

A function $g:\sE \rightarrow \R$ is said to be \emph{universally measurable} if for any $\mu \in \P(\sE)$, there is a Borel measurable function $g_{\mu}:\sE \rightarrow \R$ such that $g = g_{\mu}$ $\mu$ almost everywhere. It is said to be \emph{lower semi-analytic} if the set $\{e: g(e) < c\}$ is analytic for any $c \in \R$.  Any Borel measurable function is lower semi-analytic and any lower semi-analytic function is universally measurable. The latter property implies that the integral of any lower semi-analytic function with respect to any probability measure is well defined. We let $B^l(\sE)$ and $B^l_w(\sE)$ denote the set of all bounded lower semi-analytic functions and lower semi-analytic functions with finite $w$-norm, respectively. Since any pointwise limit of a sequence of lower semi-analytic functions is lower semi-analytic (see \citet[Theorem 1, p. 512]{Kur66}), $(B^l(\sE),\|\,\cdot\,\|)$ and $(B^l_w(\sE),\|\,\cdot\,\|_w)$ are Banach spaces.

We now state the results that will be used in the sequel.

\begin{proposition}(\citet[Proposition 7.47, p. 179]{BeSh78})\label{analytic1}
Suppose $\sE_1$ and $\sE_2$ are Borel spaces. Let $g: \sE_1 \times \sE_2 \rightarrow \R$ be lower semi-analytic. Then, $g^*(e_1) \coloneqq \inf_{e_2 \in \sE_2} g(e_1,e_2)$ is also lower semi-analytic.
\end{proposition}

\begin{proposition}(\citet[Proposition 7.48, p. 180]{BeSh78})\label{analytic2}
Suppose $\sE_1$ and $\sE_2$ as in Proposition~\ref{analytic1}. Let $g: \sE_1 \times \sE_2 \rightarrow \R$ be lower semi-analytic and $q(de_2|e_1)$ be a stochastic kernel on $\sE_2$ given $\sE_1$. Then, the function
\begin{align}
h(e_1) \coloneqq \int_{\sE_2} g(e_2) q(de_2|e_1). \nonumber
\end{align}
is lower semi-analytic.
\end{proposition}

\section{Finite State Approximations of MDPs with Compact State Space.} \label{compact}

In this section we consider \textbf{(Q1)} for the MDPs with compact state space. To distinguish compact-state MDPs from non-compact ones, the state space of the compact-state MDPs will be denoted by $\sZ$ instead of $\sX$. We impose the assumptions below on the components of the Markov decision process; additional new assumptions will be made for the average cost problem in Section~\ref{compact:sec2sub2}.

\begin{assumption}
\label{compact:as1}
\begin{itemize}
\item [  ]
\item [(a)] The one-stage cost function $c$ is in $C_b(\sZ \times \sA)$.
\item [(b)] The stochastic kernel $p(\,\cdot\,|z,a)$ is weakly continuous in $(z,a)$, i.e., for all $z$ and $a$, $p(\,\cdot\,|z_k,a_k)\rightarrow p(\,\cdot\,|z,a)$ weakly when $(z_k,a_k) \rightarrow (z,a)$.
\item [(c)] $\sZ$ and $\sA$ are compact.
\end{itemize}
\end{assumption}

Before proceeding with the main results, we first describe the procedure used to obtain finite-state models. Let $d_{\sZ}$ denote the metric on $\sZ$. Since the state space $\sZ$ is assumed to be compact and thus totally bounded, one can find a sequence $\bigl(\{z_{n,i}\}_{i=1}^{k_n}\bigr)_{n\geq1}$ of finite grids in $\sZ$ such that for all $n$,
\begin{align}
\min_{i\in\{1,\ldots,k_n\}} d_{\sZ}(z,z_{n,i}) < 1/n \text{ for all } z \in \sZ. \nonumber
\end{align}
The finite grid $\{z_{n,i}\}_{i=1}^{k_n}$ is called an $1/n$-net in $\sZ$. Let $\sZ_n \coloneqq \{z_{n,1},\ldots,z_{n,k_n}\}$ and define function $Q_n$ mapping $\sZ$ to $\sZ_n$ by
\begin{align}
Q_n(z) \coloneqq \argmin_{z_{n,i} \in \sZ_n} d_{\sZ}(z,z_{n,i}),\nonumber
\end{align}
where ties are broken so that $Q_n$ is measurable. In the literature, $Q_n$ is often called a nearest neighborhood quantizer with respect to distortion measure $d_{\sZ}$; see \citet{GrNe98}. For each $n$, $Q_n$ induces a partition $\{\S_{n,i}\}_{i=1}^{k_n}$ of the state space $\sZ$ given by
\begin{align}
\S_{n,i} = \{z \in \sZ: Q_n(z)=z_{n,i}\}, \nonumber
\end{align}
with diameter $\diam(\S_{n,i}) \coloneqq \sup_{z,y\in\S_{n,i}} d_{\sZ}(z,y) < 2/n$. Let $\{\nu_n\}$ be a sequence of
probability measures on $\sZ$ satisfying
\begin{align}
\nu_n(\S_{n,i}) > 0 \text{  for all  } i,n.  \label{compact:numeas}
\end{align}
We let $\nu_{n,i}$ be the restriction of $\nu_n$ to $\S_{n,i}$ defined by
\begin{align}
\nu_{n,i}(\,\cdot\,) \coloneqq \frac{\nu_n(\,\cdot\,)}{\nu_n(\S_{n,i})}. \nonumber
\end{align}
The measures $\nu_{n,i}$ will be used to define a sequence of finite-state MDPs, denoted as MDP$_{n}$ ($n\geq1$), to approximate the original model. To this end, for each $n$ define the one-stage cost function $c_n: \sZ_n\times\sA \rightarrow \R$ and the transition probability $p_n$ on $\sZ_n$ given $\sZ_n\times\sA$ by
\begin{align}
c_n(z_{n,i},a) &\coloneqq \int_{\S_{n,i}} c(z,a) \nu_{n,i}(dz), \nonumber \\
p_n(\,\cdot\,|z_{n,i},a) &\coloneqq \int_{\S_{n,i}} Q_n \ast p(\,\cdot\,|z,a) \nu_{n,i}(dz), \nonumber
\end{align}
where $Q_n\ast p(\,\cdot\,|z,a) \in \P(\sZ_n)$ is the pushforward of the measure $p(\,\cdot\,|z,a)$ with respect to $Q_n$; that is,
\begin{align}
Q_n\ast p(z_{n,j}|z,a) = p\bigl(\S_{n,j}|z,a\bigr), \nonumber
\end{align}
for all $z_{n,j} \in \sZ_n$. For each $n$, we define MDP$_n$ as a Markov decision process with the following components: $\sZ_n$ is the state space, $\sA$ is the action space, $p_n$ is the transition probability and $c_n$ is the one-stage cost function. History spaces, policies and cost functions are defined in a similar way as in the original model.

\subsection{Discounted Cost.}\label{compact:sec2sub1}

Here we consider \textbf{(Q1)} for the discounted cost criterion with a discount factor $\beta \in (0,1)$. Throughout this section, it is assumed that Assumption~\ref{compact:as1} holds.

Define the operator $T$ on $B(\sZ)$ by
\begin{align}
T u(z) \coloneqq \min_{a \in \sA} \biggl[ c(z,a) + \beta \int_{\sZ} u(y) p(dy|z,a) \biggr] \label{aux6}.
\end{align}
In the literature $T$ is called the \emph{Bellman optimality operator}. It can be proved that under Assumption~\ref{compact:as1}-(a)(b), $T$ is a contraction operator with modulus $\beta$ mapping $C_b(\sZ)$ into itself (see \citet[Theorem 2.8, p. 23]{Her89}); that is, $Tu \in C_b(\sZ)$ for all $u \in C_b(\sZ)$ and
\begin{align}
\|T u -  T v\| \leq \beta \|u - v\| \text{ for all } u,v \in C_b(\sZ). \nonumber
\end{align}
The following theorem is a widely known result in the theory of Markov decision processes (see again \citet[Theorem 2.8, p. 23]{Her89}) which also holds without a compactness assumption on the state space.

\begin{theorem}
\label{compact:thm1}
The value function $J^{*}$ is the unique fixed point in $C_b(\sZ)$ of the contraction operator $T$, i.e.,
\begin{align}
J^{*} = T J^{*}. \nonumber 
\end{align}
Furthermore, a deterministic stationary policy $f^{*}$ is optimal if and only if it satisfies the optimality equation, i.e.,
\begin{align}
J^{*}(z) = c(z,f^{*}(z)) + \beta \int_{\sZ} J^{*}(y) p(dy|z,f^{*}(z)). \label{compact:neq11}
\end{align}
Finally, there exists a deterministic stationary policy $f^{*}$ which is optimal, so it satisfies~(\ref{compact:neq11}).
\end{theorem}

Define, for all $n\geq1$, the operator $T_n$, which is the Bellman optimality operator for MDP$_n$, by
\begin{align}
T_n u(z_{n,i}) &\coloneqq \min_{a \in \sA} \biggl[ c_n(z_{n,i},a) + \beta \sum_{j=1}^{k_n} u(z_{n,j}) p_n(z_{n,j}|z_{n,i},a)\biggr], \nonumber
\intertext{or equivalently,}
T_n u(z_{n,i}) &= \min_{a \in \sA} \int_{\S_{n,i}} \biggl[c(z,a) + \beta \int_{\sZ} \hat{u}(y) p(dy|z,a)\biggr] \nu_{n,i}(dz), \nonumber
\end{align}
where $u: \sZ_n \rightarrow \R$ and $\hat{u}$ is the piecewise constant extension of $u$ to $\sZ$ given by $\hat{u}(z) = u\circ Q_n(z)$. For each $n$, under Assumption~\ref{compact:as1}, \citet[Theorem 2.8, p. 23]{Her89} implies the following: (i) $T_n$ is a contraction operator with modulus $\beta$ mapping $B(\sZ_n)$ $\bigl(=C_b(\sZ_n)\bigr)$ into itself, (ii) the fixed point of $T_n$ is the value function $J_n^{*}$ of MDP$_n$, and (iii) there exists an optimal stationary policy $f_n^{*}$ for MDP$_n$, which therefore satisfies the optimality equation. Hence, we have
\begin{align}
J_n^{*} = T_n J_n^{*} = T_n J_n(f_n^{*},\,\cdot\,) = J_n(f_n^{*},\,\cdot\,), \nonumber
\end{align}
where $J_n$ denotes the discounted cost for MDP$_n$. Let us extend the optimal policy $f_n^{*}$ for MDP$_n$ to $\sX$ by letting $\hf_n(z) = f_n^* \circ Q_n(z) \in \rF$.

The following theorem is the main result of this section. It states that the cost function of the policy $\hf_n$ converges to the value function $J^{*}$ as $n\rightarrow\infty$.

\begin{theorem}\label{compact:mainthm1}
The discounted cost of the policy $\hf_n$, obtained by extending the optimal policy $f_n^{*}$ of MDP$_n$ to $\sZ$,
converges to the optimal value function $J^{*}$ of the original MDP
\begin{align}
\lim_{n\rightarrow\infty} \| J(\hf_n,\,\cdot\,) - J^* \| = 0. \nonumber
\end{align}
Hence, to find a near optimal policy for the original MDP, it is sufficient to compute the optimal policy of MDP$_n$ for sufficiently large $n$, and then extend this policy to the original state space.
\end{theorem}

To prove Theorem~\ref{compact:mainthm1} we need a series of technical results.
We first define an operator $\hat{T}_n$ on $B^l(\sZ)$ by extending $T_n$ to $B^l(\sZ)$:
\begin{align}
\hat{T}_n u(z)
\coloneqq \inf_{a\in\sA} \int_{\S_{n,i_n(z)}} \biggl[c(x,a) + \beta \int_{\sZ} u(y) p(dy|x,a)\biggr] \nu_{n,i_n(z)}(dx)
, \label{compact:eq12}
\end{align}
where $i_n: \sZ \rightarrow \{1,\ldots,k_n\}$ maps $z$ to the index of the partition $\{\S_{n,i}\}$ it belongs to.
To see that this operator is well defined, let the stochastic kernel $r_n(dx|z)$ on $\sZ$ given $\sZ$ be defined as
\begin{align}
r_n(dx|z) \coloneqq \sum_{i=1}^{k_n} \nu_{n,i}(dx) 1_{\S_{n,i}}(z), \nonumber
\end{align}
where $1_B$ denotes the indicator function of the set $B$. Then, we can write the right hand side of (\ref{compact:eq12}) as
\begin{align}
\inf_{a\in\sA} \int_{\sZ} \biggl[c(x,a) + \beta \int_{\sZ} u(y) p(dy|x,a)\biggr] r_n(dx|z). \nonumber
\end{align}
Therefore, by Propositions~\ref{analytic1} and \ref{analytic2}, we can conclude that $\hat{T}_n$ maps $B^l(\sZ)$ into $B^l(\sZ)$.
Furthermore, it is a contraction operator with modulus $\beta$ which can be shown using \citet[Proposition A.2, p. 122]{Her89}. Hence, it has a unique fixed point $\hat{J}^{*}_n$ that belongs to $B(\sZ)$, and this fixed point must be constant over the sets $\S_{n,i}$ because of the averaging operation on each $\S_{n,i}$. Furthermore, since $\hat{T}_n (u \circ Q_n) = (T_n u) \circ Q_n$ for all $u \in B(\sZ_n)$, we have
\begin{align}
\hat{T}_n (J_n^{*} \circ Q_n) = (T_n J_n^{*}) \circ Q_n = J_n^{*} \circ Q_n. \nonumber
\end{align}
Hence, the fixed point of $\hat{T}_n$ is the piecewise constant extension of the fixed point of $T_n$, i.e.,
\begin{align}
\hJ^{*}_n = J_n^{*} \circ Q_n. \nonumber
\end{align}

\begin{remark}
In the rest of this paper, when we take the integral of any function with respect to $\nu_{n,i_n(z)}$, it is tacitly assumed that the integral is taken over all set $\S_{n,i_n(z)}$. Hence, we can drop $\S_{n,i_n(z)}$ in the integral for the ease of notation.
\end{remark}

We now define another operator $F_n$ on $B^l(\sZ)$ by simply interchanging the order of the infimum and the integral in (\ref{compact:eq12}), i.e.,
\begin{align}
F_n u(z) &\coloneqq \int \inf_{a \in \sA}  \biggl[c(x,a) + \beta \int_{\sZ} u(y) p(dy|x,a)\biggr] \nu_{n,i_n(z)}(dx) \nonumber \\
&= \Gamma_n T u(z), \nonumber
\end{align}
where
\begin{align}
\Gamma_n u(z) \coloneqq \int u(x) \nu_{n,i_n(z)}(dx). \nonumber
\end{align}
We note that $F_n$ is the extension (to infinite state spaces) of the operator defined in \citet[p. 236]{Roy06} for the proposed approximate value iteration algorithm. However, unlike in \citet{Roy06}, $F_n$ will serve here as an intermediate point between $T$ and $\hat{T}_n$ (or $T_n$) to solve \textbf{(Q1)} for the discounted cost. To this end, we first note that $F_n$ is a contraction operator on $B^l(\sZ)$ with modulus $\beta$. Indeed it is clear that $F_n$ maps $B^l(\sZ)$ into itself by Propositions~\ref{analytic1} and \ref{analytic2}. Furthermore, for any $u,v \in B^l(\sZ)$, we clearly have $\| \Gamma_n u - \Gamma_n v \| \leq \| u - v \|$. Hence, since $T$ is a contraction operator on $B^l(\sZ)$ with modulus $\beta$, $F_n$ is also a contraction operator on $B^l(\sZ)$ with modulus $\beta$.

\begin{remark}
Since we only assume that the stochastic kernel $p$ is weakly continuous, it is not true that $\hat{T}_n$ and $F_n$ map $B(\sZ)$ into itself (see \citet[Proposition D.5, p. 182]{HeLa96}). This is the point where we need to enlarge the set of functions on which these operators act.
\end{remark}

The following theorem states that the fixed point, say $u_n^{*}$, of $F_n$ converges to the fixed point $J^{*}$ (i.e., the value function) of $T$ as $n$ goes to infinity. Note that although $T$ is originally defined on $C_b(\sZ)$, it can be proved that $T$, when acting on $B^l(\sZ)$, maps $B^l(\sZ)$ into itself.

\begin{theorem}\label{compact:thm2}
If $u_n^{*}$ is the unique fixed point of $F_n$, then $\lim_{n\rightarrow\infty} \| u_n^{*} - J^{*} \| = 0$.
\end{theorem}

The proof of Theorem~\ref{compact:thm2} requires two lemmas.

\begin{lemma}\label{compact:lemma1}
For any $u \in B^l(\sZ)$, we have
\begin{align}
\| u - \Gamma_n u \| \leq 2 \inf_{r \in \sZ^{k_n}} \| u - \Phi_{r} \|, \nonumber
\end{align}
where $\Phi_{r}(z) = \Sigma_{i=1}^{k_n} r_i 1_{S_{n,i}}(z)$, $r = (r_1,\cdots,r_{k_n})$.
\end{lemma}

\proof{Proof.}
Fix any $r \in \sZ^{k_n}$. Then, using the identity $\Gamma_n \Phi_{r} = \Phi_{r}$, we obtain
\begin{align}
\| u - \Gamma_n u \| &\leq \| u - \Phi_{r} \| + \| \Phi_{r} - \Gamma_n u \| \nonumber \\
&= \| u - \Phi_{r} \| + \| \Gamma_n \Phi_{r} - \Gamma_n u \| \nonumber \\
&\leq \| u - \Phi_{r} \| + \| \Phi_{r} - u \|. \nonumber
\end{align}
Since $r$ is arbitrary, this completes the proof.\Halmos
\endproof

Notice that because of the operator $\Gamma_n$, the fixed point $u_n^{*}$ of $F_n$ must be constant over the sets $\S_{n,i}$. We use this property to prove the next lemma.

\begin{lemma}\label{compact:nlemma2}
We have
\begin{align}
\| u_n^{*} - J^{*} \| \leq \frac{2}{1-\beta} \inf_{r \in \sZ^{k_n}} \| J^{*} - \Phi_{r} \|. \nonumber
\end{align}
\end{lemma}

\proof{Proof.}
Note that $\Gamma_n u_n^{*} = u_n^{*}$ since $u_n^{*}$ is constant over the sets $\S_{n,i}$. Then, we have
\begin{align}
\| u_n^{*} - J^* \| &\leq \| u_n^{*} - \Gamma_n J^* \| + \| \Gamma_n J^* - J^* \| \nonumber \\
&= \| F_n u_n^{*} - \Gamma_n T J^* \| + \| \Gamma_n J^* - J^* \| \nonumber \\
&= \| \Gamma_n T u_n^{*} - \Gamma_n T J^* \| + \| \Gamma_n J^* - J^* \| \text{  (by the definition of $F_n$)} \nonumber \\
&\leq \| T u_n^{*} - T J^* \| + \| \Gamma_n J^* - J^* \|  \text{ (since $\|\Gamma_n u - \Gamma_n v\| \leq \| u - v\|$)}\nonumber \\
&\leq \beta \| u_n^{*} - J^* \| + \| \Gamma_n J^* - J^* \|. \nonumber
\end{align}
Hence, we obtain $\| u_n^{*} - J^* \| \leq \frac{1}{1-\beta} \| \Gamma_n J^* - J^* \|$. The result now follows from Lemma~\ref{compact:lemma1}. \Halmos
\endproof

\proof{Proof of Theorem~\ref{compact:thm2}.}
Recall that since $\sZ$ is compact, the function $J^*$ is uniformly continuous and $\diam(\S_{n,i}) < 2/n$ for all $i=1,\ldots,k_n$. Hence, $\lim_{n\rightarrow\infty} \inf_{r \in \sZ^{k_n}} \| J^* - \Phi_{r} \| = 0$ which completes the proof in view of Lemma~\ref{compact:nlemma2}. \Halmos
\endproof

The next step is to show that the fixed point $\hJ^*_n$ of $\hat{T}_n$ converges to the fixed point $J^*$ of $T$. To this end, we first prove the following result.

\begin{lemma}\label{compact:nlemma3}
For any $u \in C_b(\sZ)$, $\| \hat{T}_n u - F_n u \| \rightarrow 0$ as $n\rightarrow\infty$.
\end{lemma}

\proof{Proof.}
Note that since $\int_{\sZ} u(x) p(dx|y,a)$ is continuous as a function of $(y,a)$ by Assumption~\ref{compact:as1}-(b), it is sufficient to prove that for any $l \in C_b(\sZ\times\sA)$
\begin{align}
&\biggl \| \min_a \int l(y,a) \nu_{n,i_n(z)}(dy) - \int \min_a l(y,a) \nu_{n,i_n(z)}(dy) \biggr \| \nonumber \\
&\phantom{xxxxxxxx}\coloneqq \sup_{z\in\sZ} \phantom{i} \biggl|\min_a \int l(y,a) \nu_{n,i_n(z)}(dy) - \int \min_a l(y,a) \nu_{n,i_n(z)}(dy) \biggr |\rightarrow 0 \nonumber
\end{align}
as $n\rightarrow\infty$. Fix any $\varepsilon > 0$. Define $\{z_i\}_{i=1}^\infty \coloneqq \bigcup_n \sZ_n$ and let $\{a_i\}_{i=1}^\infty$ be a sequence in $\sA$ such that $\min_{a\in\sA} l(z_i,a) = l(z_i,a_i)$; such $a_i$ exists for each $z_i$ because $l(z_i,\,\cdot\,)$ is continuous and $\sA$ is compact. Define $g(y) \coloneqq \min_{a\in\sA} l(y,a)$, which can be proved to be continuous, and therefore uniformly continuous since $\sZ$ is compact. Thus by the uniform continuity of $l$, there exists $\delta > 0$ such that $d_{\sZ\times\sA}\bigl((y,a),(y',a')\bigr) < \delta$ implies $|g(y) - g(y')| < \varepsilon/2$ and $|l(y,a) - l(y',a')| < \varepsilon/2$. Choose $n_0$ such that $2/n_0 < \delta$. Then for all $n\geq n_0$, $\max_{i\in\{1,\ldots,k_n\}} \diam(\S_{n,i}) < 2/n < \delta$. Hence, for all $y \in \S_{n,i}$ we have $|l(y,a_i) - \min_{a\in\sA} l(y,a)| \leq |l(y,a_i) - l(z_i,a_i)| + |\min_{a\in\sA} l(z_i,a) - \min_{a\in\sA} l(y,a)| = |l(y,a_i) - l(z_i,a_i)| + |g(z_i) - g(y)| < \varepsilon$. This implies
\begin{align}
&\biggl \| \min_a \int l(y,a) \nu_{n,i_n(z)}(dy) - \int \min_a l(y,a) \nu_{n,i_n(z)}(dy) \biggr \| \nonumber \\
&\phantom{xxxxxxxxxx}\leq \biggl \| \int l(y,a_i) \nu_{n,i_n(z)}(dy) - \int \min_a l(y,a) \nu_{n,i_n(z)}(dy) \biggr \| \nonumber \\
&\phantom{xxxxxxxxxx}\leq \sup_{z \in \sZ} \int \sup_{y \in \S_{n,i_n(z)}} \bigl|l(y,a_i) - \min_a l(y,a) \bigr| \nu_{n,i_n(z)}(dy)  < \varepsilon. \nonumber
\end{align}
This completes the proof. \Halmos
\endproof

\begin{theorem}\label{compact:thm3}
The fixed point $\hat{J}_n^{*}$ of $\hat{T}_n$ converges to the fixed point $J^{*}$ of $T$.
\end{theorem}

\proof{Proof.}
We have
\begin{align}
\| \hJ_n^* - J^* \| &\leq \| \hat{T}_n \hJ_n^* - \hat{T}_n J^* \| + \| \hat{T}_n J^* - F_n J^* \|
            + \| F_n J^* - F_n u_n^{*} \| \nonumber \\
&\phantom{xxxxxxxxxxxxxxxxxxxxxxxxxxxxxxxx}+ \|F_n u_n^{*} - J^*\| \nonumber \\
&\leq \beta \| \hJ_n^* -J^* \| + \| \hat{T}_n J^* - F_n J^* \|
+\beta \| J^* - u_n^{*} \| + \| u_n^{*} - J^* \|. \nonumber
\end{align}
Hence
\begin{align}
\| \hJ_n^* - J^* \| \leq \frac{\| \hat{T}_n J^* - F_n J^* \| + (1+\beta) \| J^* - u_n^{*} \|}{1-\beta}. \nonumber
\end{align}
The theorem now follows from Theorem~\ref{compact:thm2} and Lemma~\ref{compact:nlemma3}. \Halmos
\endproof

Recall the optimal stationary policy $f_n^*$ for MDP$_n$ and its extension $\hf_n(z) = f_n^*\circ Q_n(z)$ to $\sZ$. Since $\hJ_n^{*} = J_n^{*} \circ Q_n$, it is straightforward to prove that $\hf_n$ is the optimal selector of $\hat{T}_n \hJ_n^*$; that is,
\begin{align}
\hat{T}_n \hJ_n^* = \hJ_n^* = \hat{T}_{\hf_n} \hJ_n^*, \nonumber
\end{align}
where $\hat{T}_{\hf_n}$ is defined as
\begin{align}
\hat{T}_{\hf_n} u(z) &\coloneqq \int \biggl[c(x,\hf_n(x)) + \beta \int_{\sZ} u(y) p(dy|x,\hf_n(x))\biggr] \nu_{n,i_n(z)}(dx). \nonumber
\intertext{Define analogously}
T_{\hf_n} u(z) &\coloneqq c(z,\hf_n(z)) + \beta \int_{\sZ} u(y) p(dy|z,\hf_n(z)). \nonumber
\end{align}
It can be proved that both $\hat{T}_{\hf_n}$ and $T_{\hf_n}$ are contraction operators on $B^l(\sZ)$ with modulus $\beta$, and it is known that the fixed point of $T_{\hf_n}$ is the true cost function of the stationary policy $\hf_n$ (i.e., $J(\hf_n,z)$).

\begin{lemma}\label{compact:lemma2}
$\| \hat{T}_{\hf_n} u - T_{\hf_n} u \| \rightarrow 0$ as $n \rightarrow \infty$, for any $u \in C_b(\sZ)$.
\end{lemma}

\proof{Proof.}
The statement follows from the uniform continuity of the function $c(z,a) + \beta \int_{\sZ} u(y) p(dy|z,a)$ and the fact that $\hf_n$ is constant over the sets $\S_{n,i}$. \Halmos
\endproof

Now, we prove the main result of this section.

\proof{Proof of Theorem~\ref{compact:mainthm1}.}
We have
\begin{align}
\| J(\hf_n,\,\cdot\,) - J^* \| &\leq \| T_{\hf_n} J(\hf_n,\,\cdot\,) - T_{\hf_n} J^* \| + \| T_{\hf_n} J^* - \hat{T}_{\hf_n} J^* \| + \| \hat{T}_{\hf_n} J^* - \hat{T}_{\hf_n} \hJ^*_n \| + \| \hJ^*_n - J^* \| \nonumber \\
& \leq \beta \| J(\hf_n,\,\cdot\,) - J^* \| + \| T_{\hf_n} J^* - \hat{T}_{\hf_n} J^* \| + \beta \| J^* - \hJ^*_n \| + \| \hJ^*_n - J^* \|. \nonumber
\end{align}
Hence, we obtain
\begin{align}
\hspace{-2pt}\| J(\hf_n,\,\cdot\,) - J^* \| \leq \frac{\| T_{\hf_n} J^* - \hat{T}_{\hf_n} J^* \| + (1+\beta) \| \hJ^*_n - J^* \|}{1-\beta}. \nonumber
\end{align}
The result follows from Lemma~\ref{compact:lemma2} and Theorem~\ref{compact:thm3}. \Halmos
\endproof

\subsection{Average Cost.}\label{compact:sec2sub2}

In this section we impose some new conditions on the components of the original MDP in addition to Assumption~\ref{compact:as1} to solve \textbf{(Q1)} for the average cost. A version of the first two conditions was imposed in \citet{Veg03,JaNo06} to show the existence of the solution to the Average Cost Optimality Equation (ACOE) and the optimal stationary policy.

\begin{assumption}
\label{compact:as2}
Suppose Assumption~\ref{compact:as1} holds with item (b) replaced by condition (f) below. In addition, there exist a non-trivial finite measure $\zeta$ on $\sZ$, a nonnegative measurable function $\theta$ on $\sZ \times \sA$, and a constant $\lambda \in (0,1)$ such that for all $(z,a) \in \sZ\times\sA$
\begin{itemize}
\item [(d)] $p(B|z,a) \geq \zeta(B) \theta(z,a)$ for all $B \in \B(\sZ)$,
\item [(e)] $\frac{1-\lambda}{\zeta(\sZ)} \leq  \theta(z,a)$,
\item [(f)] The stochastic kernel $p(\,\cdot\,|z,a)$ is continuous in $(z,a)$ with respect to the total variation distance.
\end{itemize}
\end{assumption}

Throughout this section, it is assumed that Assumption~\ref{compact:as2} holds.
Observe that any deterministic stationary policy $f$ defines a stochastic kernel $p(\,\cdot\,|z,f(z))$ on $\sZ$ given $\sZ$ which is the transition probability of the Markov chain $\{z_t\}_{t=1}^{\infty}$ (state process) induced by $f$. For any $t\geq1$, let us write $p^t(\,\cdot\,|z,f(z))$ to denote the $t$-step transition probability of this Markov chain given the initial point $z$; that is, $p^t(\,\cdot\,|z,f(z))$ is recursively defined as
\begin{align}
p^{t+1}(\,\cdot\,|z,f(z)) = \int_{\sZ} p(\,\cdot\,|x,f(x)) p^t(dx|z,f(z)). \nonumber
\end{align}

To study average cost optimal control problems, it is in general assumed that there exists an invariant distribution under any stationary control policy, so that the average cost of any stationary policy can be written as an integral of the one-stage cost function with respect to this invariant distribution. With this representation, one can then deduce the optimality of stationary policies using the linear programming or the convex analytic methods (see \citet{HeLa96,Bor02}). However, to solve the approximation problem for the average cost, we need, in addition to the existence of an invariant distribution, the convergence of $t$-step transition probabilities to the invariant distribution, at some rate, for both the original and the reduced problems. Therefore, it is crucial to impose proper conditions on the original model so that, on the one hand, they guarantee the convergence of $t$-step transition probabilities to the invariant distribution for all stationary policies for the original system and, on the other hand, one is able to show that similar conditions are satisfied by the reduced problems. Conditions (d) and (e) in Assumption~\ref{compact:as2} are examples of such conditions which were also used in the literature extensively.  
Indeed, if we define the weight function $w \equiv 1$, then condition (e) corresponds to the so-called `drift inequality': for all $(z,a) \in \sZ \times \sA$ \begin{align}
\int_{\sZ} w(y) p(dy|z,a) &\leq \lambda w(z) + \zeta(w) \theta(z,a), \nonumber
\end{align}
and condition (d) corresponds to the so-called `minorization' condition, both of which were used in literature for studying geometric ergodicity of Markov chains (see \citet{HeLa99,MeTw93}, and references therein).

The following theorem is a consequence of \citet[Theorem 3.3]{Veg03}, \citet[Lemma 3.4]{GoHe95}, and \citet[Theorem 3]{JaNo06}, which also holds with Assumption~\ref{compact:as2}-(f) replaced by Assumption~\ref{compact:as1}-(b).

\begin{theorem}\label{compact:thm4}
For any $f \in \rF$, the stochastic kernel $p(\,\cdot\,|z,f(z))$ is positive Harris recurrent with unique invariant probability measure $\mu_f$. Therefore, we have
\begin{align}
V(f,z) = \int_{\sZ} c(z,f(z)) \mu_f(dz) \eqqcolon \rho_f. \nonumber
\end{align}
The Markov chain $\{z_t\}_{t=1}^{\infty}$ induced by $f$ is geometrically ergodic; that is, there exist positive real numbers $R$ and $\kappa < 1$ such that for every $z \in \sZ$
\begin{align}
\sup_{f \in \rF} \| p^t(\,\cdot\,|z,f(z)) - \mu_f \|_{TV} \leq  R \kappa^t, \nonumber
\end{align}
where $R$ and $\kappa$ continuously depend on $\zeta(\sZ)$ and $\lambda$.
Finally, there exist $f^{*} \in \rF$ and  $h^{*} \in B(\sZ)$ such that the triplet $(h^{*},f^{*},\rho_{f^{*}})$ satisfies the average cost optimality equality (ACOE), i.e.,
\begin{align}
\rho_{f^{*}} + h^{*}(z) &= \min_{a\in\sA} \biggl[ c(z,a) + \int_{\sZ} h^{*}(y) p(dy|z,a) \biggr] \nonumber \\
&=  c(z,f^{*}(z)) + \int_{\sZ} h^{*}(y) p(dy|z,f^{*}(z)), \nonumber
\end{align}
and therefore,
\begin{align}
\inf_{\pi \in \Pi} V(\pi,z) \eqqcolon V^{*}(z) = \rho_{f^{*}}. \nonumber
\end{align}
\end{theorem}

For each $n$, define the one-stage cost function $b_n: \sZ\times\sA \rightarrow [0,\infty)$ and the stochastic kernel $q_n$ on $\sZ$ given $\sZ\times\sA$ as
\begin{align}
b_n(z,a) &\coloneqq \int c(x,a) \nu_{n,i_n(z)}(dx), \nonumber \\
q_n(\,\cdot\,|z,a) &\coloneqq \int p(\,\cdot\,|x,a) \nu_{n,i_n(z)}(dx). \nonumber
\end{align}
Observe that $c_n$ (i.e., the one stage cost function of MDP$_n$) is the restriction of $b_n$ to $\sZ_n$, and $p_n$ (i.e., the stochastic kernel of MDP$_n$) is the pushforward of the measure $q_n$ with respect to $Q_n$; that is, $c_n(z_{n,i},a)=b_n(z_{n,i},a)$ for all $i=1,\ldots,k_n$ and $p_n(\,\cdot\,|z_{n,i},a) = Q_n \ast q_n(\,\cdot\,|z_{n,i},a)$.

For each $n$, let $\widehat{\text{MDP}}_n$ be defined as a Markov decision process with the following components: $\sZ$ is the state space, $\sA$ is the action space, $q_n$ is the transition probability, and $c$ is the one-stage cost function. Similarly, let $\widetilde{\text{MDP}}_n$ be defined as a Markov decision process with the following components: $\sZ$ is the state space, $\sA$ is the action space, $q_n$ is the transition probability, and $b_n$ is the one-stage cost function. History spaces, policies and cost functions are defined in a similar way as before. The models $\widehat{\text{MDP}}_n$ and $\widetilde{\text{MDP}}_n$ are used as  transitions between the original MDP and MDP$_n$ in a similar way as the operators $F_n$ and $\hat{T}_n$ were used as transitions between $T$ and $T_n$ for the discounted cost. We note that a similar technique was used in the proof of \citet[Theorem 2]{Ort07}, which studied the approximation problem for finite state-action MDPs. In \citet{Ort07} the one-stage cost function is first perturbed and then the transition probability is perturbed. We first perturb the transition probability and then the cost function. However, our proof method is otherwise quite different from that of \citet[Theorem 2]{Ort07} since \citet{Ort07} assumes finite state and action spaces.

We note that a careful analysis of $\widetilde{\text{MDP}}_n$ reveals that its Bellman optimality operator is essentially the operator $\hat{T}_n$. Hence, the value function of $\widetilde{\text{MDP}}_n$ is the piecewise constant extension of the value function of MDP$_n$ for the discounted cost. A similar conclusion will be made for the average cost in Lemma~\ref{compact:lemma0}.

First, notice that if we define
\begin{align}
\theta_n(z,a) &\coloneqq \int \theta(y,a) \nu_{n,i_n(z)}(dy), \nonumber \\
\zeta_n &\coloneqq Q_n \ast \zeta  \text{ (i.e., pushforward of $\zeta$ with respect to $Q_n$)}, \nonumber
\end{align}
then it is straightforward to prove that for all $n$, both $\widehat{\text{MDP}}_n$ and $\widetilde{\text{MDP}}_n$ satisfy Assumption~\ref{compact:as2}-(d),(e) when $\theta$ is replaced by $\theta_n$, and Assumption~\ref{compact:as2}-(d),(e) is true for MDP$_n$ when $\theta$ and $\zeta$ are replaced by the restriction of $\theta_n$ to $\sZ_n$ and $\zeta_n$, respectively.

Hence, Theorem~\ref{compact:thm4} holds (with the same $R$ and $\kappa$) for $\widehat{\text{MDP}}_n$, $\widetilde{\text{MDP}}_n$, and MDP$_n$ for all $n$. Therefore, we denote by $\hat{f}_n^{*}$, $\tilde{f}_n^{*}$ and $f_n^{*}$ the optimal stationary policies of $\widehat{\text{MDP}}_n$, $\widetilde{\text{MDP}}_n$, and MDP$_n$ with the corresponding average costs $\hat{\rho}^n_{\hat{f}_n^{*}}$, $\tilde{\rho}^n_{\tilde{f}_n^{*}}$ and $\rho^n_{f_n^{*}}$, respectively.

Furthermore, we also write $\hat{\rho}^n_f$, $\tilde{\rho}^n_f$, and $\rho^n_f$ to denote the average cost of any stationary policy $f$ for $\widehat{\text{MDP}}_n$, $\widetilde{\text{MDP}}_n$, and MDP$_n$, respectively. The corresponding invariant probability measures are also denoted in a similar manner, with $\mu$ replacing $\rho$.

The following lemma essentially says that MDP$_n$ and $\widetilde{\text{MDP}}_n$ are not very different.

\begin{lemma}\label{compact:lemma0}
The stationary policy given by the piecewise constant extension of the optimal policy $f_n^{*}$ of MDP$_n$ to $\sZ$ (i.e., $f_n^{*} \circ Q_n$) is optimal for $\widetilde{\text{MDP}}_n$ with the same cost function $\rho^n_{f_n^{*}}$. Hence, $\tilde{f}_n^{*} = f_n^{*} \circ Q_n$ and $\tilde{\rho}^n_{\tilde{f}_n^{*}} = \rho^n_{f_n^{*}}$.
\end{lemma}

\proof{Proof.}
Note that by Theorem~\ref{compact:thm4} there exists $h_n^{*} \in B(\sZ_n)$ such that the triplet $(h_n^{*},f_n^{*},\rho^n_{f_n^{*}})$ satisfies the ACOE for MDP$_n$. But it is straightforward to show that the triplet $(h_n^{*}\circ Q_n,f_n^{*}\circ Q_n,\rho^n_{f_n^{*}})$ satisfies the ACOE for $\widetilde{\text{MDP}}_n$.
By \citet[Lemma 5.2]{GoHe95}, this implies that $f_n^{*}\circ Q_n$ is an optimal stationary policy for $\widetilde{\text{MDP}}_n$ with cost function $\rho^n_{f_n^{*}}$. Hence $\tilde{f}_n^{*} = f_n^{*}\circ Q_n$ and $\tilde{\rho}^n_{\tilde{f}_n^{*}} = \rho^n_{f_n^{*}}$. \Halmos
\endproof

The following theorem is the main result of this section. It states that if one applies the piecewise constant extension of the optimal stationary policy of MDP$_n$ to the original MDP, the resulting cost function will converge to the value function of the original MDP.

\begin{theorem}\label{compact:thm6}
The average cost of the optimal policy $\tilde{f}^{*}_n$ for $\widetilde{\text{MDP}}_n$, obtained by extending the optimal policy $f_n^{*}$ of MDP$_n$ to $\sZ$,
converges to the optimal value function $J^{*}=\rho_{f^{*}}$ of the original MDP, i.e.,
\begin{align}
\lim_{n\rightarrow\infty} |\rho_{\tilde{f}_n^{*}} - \rho_{f^{*}}| = 0. \nonumber
\end{align}
Hence, to find a near optimal policy for the original MDP, it is sufficient to compute the optimal policy of MDP$_n$ for sufficiently large $n$, and then extend this policy to the original state space.
\end{theorem}

To show the statement of Theorem~\ref{compact:thm6} we will prove a series of auxiliary results.
\begin{lemma}\label{compact:prop3}
For all $t\geq1$ we have
\begin{align}
\lim_{n\rightarrow\infty} \sup_{(y,f) \in \sZ\times\rF} \bigl\|p^t(\,\cdot\,|y,f(y)) - q_n^t(\,\cdot\,|y,f(y)) \bigr\|_{TV} = 0. \nonumber
\end{align}
\end{lemma}

\proof{Proof.}
We will prove the lemma by induction. Note that if one views the stochastic kernel $p(\,\cdot\,|z,a)$ as a mapping from $\sZ\times\sA$ to $\P(\sZ)$, then Assumption~\ref{compact:as2}-(f) implies that this mapping is continuous, and therefore uniformly continuous, when $\P(\sZ)$ is equipped with the metric induced by the total variation distance.

For $t=1$ the claim holds by the following argument:
\begin{align}
\sup_{(y,f) \in \sZ\times\rF} \bigl\|p(\,\cdot\,|y,f(y)) - q_n(\,\cdot\,|y,f(y)) \bigr\|_{TV} &\coloneqq 2 \sup_{(y,f) \in \sZ\times\rF} \sup_{D\in\B(\sZ)} \bigl|p(D|y,f(y)) - q_n(D|y,f(y))\bigr| \nonumber \\
&\leq 2 \sup_{(y,f) \in \sZ\times\rF} \sup_{D\in\B(\sZ)} \int \bigl|p(D|y,f(y)) - p(D|z,f(y))\bigr| \hspace{3pt} \nu_{n,i_n(y)}(dz) \nonumber \\
&\leq \sup_{(y,f) \in \sZ\times\rF} \int \bigl\|p(\,\cdot\,|y,f(y)) - p(\,\cdot\,|z,f(y))\bigr\|_{TV} \nu_{n,i_n(y)}(dz) \nonumber \\
&\leq \sup_{y\in\sZ} \sup_{(z,a) \in \S_{n,i_n(y)}\times\sA} \bigl\|p(\,\cdot\,|y,a) - p(\,\cdot\,|z,a)\bigr\|_{TV}. \nonumber
\end{align}
As the mapping $p(\,\cdot\,|z,a):\sZ\times\sA \rightarrow \P(\sZ)$ is uniformly continuous with respect to the total variation distance and $\max_{n,i} \diam(\S_{n,i}) \rightarrow 0$ as $n\rightarrow\infty$, the result follows. Assume the claim is true for $t\geq1$. Then we have
\begin{align}
&\sup_{(y,f) \in \sZ\times\rF} \bigl\| p^{t+1}(\,\cdot\,|y,f(y)) - q_n^{t+1}(\,\cdot\,|y,f(y)) \bigr\|_{TV} \nonumber \\
&\phantom{xxxxx}\coloneqq \sup_{(y,f) \in \sZ\times\rF} \sup_{\|g\| \leq 1} \biggl| \int_{\sZ} g(x) p^{t+1}(dx|y,f(y)) - \int_{\sZ} g(x) q_n^{t+1}(dx|y,f(y)) \biggr| \nonumber \\
&\phantom{xxxxx}\leq \sup_{(y,f) \in \sZ\times\rF} \biggl( \sup_{\|g\|\leq 1} \biggl| \int_{\sZ} \int_{\sZ} g(x) p(dx|z,f(z)) p^t(dz|y,f(y)) - \int_{\sZ} \int_{\sZ} g(x) p(dx|z,f(z)) q_n^t(dz|y,f(y)) \biggr| \nonumber \\
&\phantom{xxxxxxxxxx} + \sup_{\|g\|\leq 1} \biggl| \int_{\sZ} \int_{\sZ} g(x) p(dx|z,f(z)) q_n^t(dz|y,f(y)) - \int_{\sZ} \int_{\sZ} g(x) q_n(dx|z,f(z)) q_n^t(dz|y,f(y)) \biggr| \biggl) \nonumber \\
&\phantom{xxxxx}\leq \hspace{-10pt} \sup_{(y,f) \in \sZ\times\rF} \bigl\| p^t(\,\cdot\,|y,f(y)) - q_n^t(\,\cdot\,|y,f(y)) \bigr\|_{TV} \hspace{-5pt}+
\hspace{-8pt} \sup_{(z,f) \in \sZ\times\rF} \bigl\| p(\,\cdot\,|z,f(z)) - q_n(\,\cdot\,|z,f(z)) \bigr\|_{TV} \label{aux11}
\end{align}
where the last inequality follows from the following property of the total variation distance: for any $h \in \B(\sZ)$ and $\mu,\nu \in \P(\sZ)$ we have $\bigl| \int_{\sZ} h(z) \mu(dz) - \int_{\sZ} h(z) \nu(dz) \bigr| \leq \|h\| \|\mu-\nu\|_{TV}$.
By the first step of the proof and the induction hypothesis, the last term converges to zero as $n\rightarrow\infty$. This completes the proof. \Halmos
\endproof

\begin{remark}\label{remarktotal}
This is the point where we need the continuity of the transition probability $p$ with respect to the total variation distance. If we assume that the stochastic kernel $p$ is only weakly or setwise continuous, then it does not seem possible to prove a result similar to Lemma~\ref{compact:prop3} for the weak and the setwise topologies.
\end{remark}

Using Lemma~\ref{compact:prop3} we prove the following result.
\begin{lemma}\label{compact:prop4}
We have $\sup_{f \in \rF} | \hat{\rho}^n_f - \rho_f | \rightarrow 0$ as $n\rightarrow\infty$, where $\hat{\rho}^n_f$ is the cost function of the policy $f$ for $\widehat{\text{MDP}}_n$ and $\rho_f$ is the cost function of the policy $f$ for the original MDP.
\end{lemma}

\proof{Proof.}
For any $t\geq1$ and $y \in \sZ$ we have
\begin{align}
\sup_{f \in \rF} |\hat{\rho}_f^n - \rho_f| &= \sup_{f \in \rF} \biggl | \int_{\sZ} c(z,f(z)) \hat{\mu}_f^n(dz) - \int_{\sZ} c(z,f(z)) \mu_f(dz) \biggr | \nonumber \\
&\leq \sup_{f \in \rF} \biggl | \int_{\sZ} c(z,f(z)) \hat{\mu}_f^n(dz) - \int_{\sZ} c(z,f(z)) q_n^t(dz|y,f(y)) \biggr | \nonumber \\
&\phantom{xxxxxxx}+ \sup_{f \in \rF} \biggl | \int_{\sZ} c(z,f(z)) q_n^t(dz|y,f(y)) - \int_{\sZ} c(z,f(z)) p^t(dz|y,f(y)) \biggr | \nonumber \\
&\phantom{xxxxxxxxxxxxxxxx}+ \sup_{f \in \rF} \biggl | \int_{\sZ} c(z,f(z)) p^t(dz|y,f(y)) - \int_{\sZ} c(z,f(z)) \mu_f(dz) \biggr | \nonumber \\
&\leq 2 R \kappa^t \|c\|  + \|c\| \sup_{(y,f)\in\sZ\times\rF} \bigl\| q_n^t(\,\cdot\,|y,f(y)) - p^t(\,\cdot\,|y,f(y)) \bigr\|_{TV} \text{ (by Theorem~\ref{compact:thm4}-(ii))}, \nonumber
\end{align}
where $R$ and $\kappa$ are the constants in Theorem~\ref{compact:thm4}. Then, the result follows from Lemma~\ref{compact:prop3}. \Halmos
\endproof

The following theorem states that the value function of $\widehat{\text{MDP}}_n$ converges to the value function of the original MDP.

\begin{lemma}\label{compact:lemma3}
We have $|\hat{\rho}^n_{\hat{f}_n^{*}}-\rho_{f^{*}}| \rightarrow 0$ as $n\rightarrow\infty$.
\end{lemma}

\proof{Proof.}
Notice that
\begin{align}
|\hat{\rho}^n_{\hat{f}_n^{*}} - \rho_{f^{*}}| &= \max (\hat{\rho}^n_{\hat{f}_n^{*}} - \rho_{f^{*}}, \rho_{f^{*}} - \hat{\rho}^n_{\hat{f}_n^{*}}) \nonumber \\
&\leq \max (\hat{\rho}^n_{f^{*}} - \rho_{f^{*}}, \rho_{\hat{f}^{*}_n} - \hat{\rho}^n_{\hat{f}_n^{*}}) \nonumber \\
&\leq \sup_f |\hat{\rho}^n_f - \rho_f|. \nonumber
\end{align}
Then, the result follows from Lemma~\ref{compact:prop4}. \Halmos
\endproof

\begin{lemma}\label{compact:prop5}
We have $\sup_{f\in\rF} | \tilde{\rho}^n_f - \hat{\rho}^n_f | \rightarrow 0$ as $n\rightarrow\infty$.
\end{lemma}

\proof{Proof.}
It is straightforward to show that $b_n \rightarrow c$ uniformly. Since the probabilistic structure of $\widetilde{\text{MDP}}_n$ and $\widehat{\text{MDP}}_n$ are the same (i.e., $\hat{\mu}^n_f = \tilde{\mu}^n_f$ for all $f$), we have
\begin{align}
\sup_{f\in\rF} | \tilde{\rho}^n_f - \hat{\rho}^n_f | &= \sup_{f\in\rF} \biggl | \int_{\sZ} b_n(z,f(z)) \hat{\mu}^n_f(dz) -  \int_{\sZ} c(z,f(z)) \hat{\mu}^n_f(dz) \biggr | \nonumber \\
&\leq \sup_{f\in\rF} \int_{\sZ} |b_n(z,f(z))-c(z,f(z))| \hat{\mu}^n_f(dz) \nonumber \\
&\leq \|b_n - c\|.\nonumber
\end{align}
This completes the proof. \Halmos
\endproof

The next lemma states that the difference between the value functions of $\widetilde{\text{MDP}}_n$ and $\widehat{\text{MDP}}_n$ converges to zero.

\begin{lemma}\label{compact:prop6}
We have $|\tilde{\rho}^n_{\tilde{f}_n^{*}}-\hat{\rho}^n_{\hat{f}_n^{*}}|\rightarrow 0$ as $n\rightarrow\infty$.
\end{lemma}

\proof{Proof.}
See the proof of Lemma~\ref{compact:lemma3}. \Halmos
\endproof

The following result states that if we apply the optimal policy of $\widetilde{\text{MDP}}_n$ to $\widehat{\text{MDP}}_n$, then the resulting cost converges to the value function of $\widehat{\text{MDP}}_n$.

\begin{lemma}\label{compact:prop7}
We have $|\hat{\rho}^n_{\tilde{f}_n^{*}} - \hat{\rho}^n_{\hat{f}_n^{*}}| \rightarrow 0$ as $n\rightarrow\infty$.
\end{lemma}

\proof{Proof.}
Since $|\hat{\rho}^n_{\tilde{f}_n^{*}} - \hat{\rho}^n_{\hat{f}_n^{*}}| \leq |\hat{\rho}^n_{\tilde{f}_n^{*}} - \tilde{\rho}^n_{\tilde{f}_n^{*}}| + |\tilde{\rho}^n_{\tilde{f}_n^{*}} - \hat{\rho}^n_{\hat{f}_n^{*}}|$, then the result follows from
Lemmas~\ref{compact:prop5} and \ref{compact:prop6}. \Halmos
\endproof

Now, we are ready to prove the main result of this section.

\proof{Proof of Theorem~\ref{compact:thm6}.}
We have $|\rho_{\tilde{f}_n^{*}} - \rho_{f^{*}}| \leq |\rho_{\tilde{f}_n^{*}} - \hat{\rho}^n_{\tilde{f}_n^{*}}| + |\hat{\rho}^n_{\tilde{f}_n^{*}} - \hat{\rho}^n_{\hat{f}_n^{*}}| + |\hat{\rho}^n_{\hat{f}_n^{*}} - \rho^n_{f^{*}}|$. The result now follows from Lemmas~\ref{compact:prop4}, \ref{compact:prop7} and \ref{compact:lemma3}. \Halmos
\endproof

\section{Finite State Approximations of MDPs with Non-Compact State Space.}\label{sec2}

In this section we consider \textbf{(Q1)} for noncompact state MDPs with unbounded one-stage cost. To solve \textbf{(Q1)}, we use the following strategy: (i) first, we define a sequence of compact-state MDPs to approximate the original MDP, (ii) we use Theorems~\ref{compact:mainthm1} and \ref{compact:thm6} to approximate the compact-state MDPs by finite-state models, and (iii) we prove the convergence of the finite-state models to the original model. In fact, steps (ii) and (iii) will be accomplished simultaneously.

We impose the assumptions below on the components of the Markov decision process; additional assumptions will be imposed for the average cost problem. With the exception of the local compactness of the state space, these are the usual assumptions used in the literature for studying Markov decision processes with unbounded cost.

\begin{assumption}
\label{as1}
\begin{itemize}
\item [ ]
\item [(a)] The one-stage cost function $c$ is continuous.
\item [(b)] The stochastic kernel $p(\,\cdot\,|x,a)$ is weakly continuous in $(x,a)$.
\item [(c)] $\sX$ is locally compact and $\sA$ is compact.
\item [(d)] There exist nonnegative real numbers $M$ and $\alpha \in [1,\frac{1}{\beta})$, and a continuous weight function $w:\sX\rightarrow[1,\infty)$ such that for each $x \in \sX$, we have
\begin{align}
\sup_{a \in \sA} |c(x,a)| &\leq M w(x), \label{eq1} \\
\sup_{a \in \sA} \int_{\sX} w(y) p(dy|x,a) &\leq \alpha w(x), \label{eq2}
\end{align}
and $\int_{\sX} w(y) p(dy|x,a)$ is continuous in $(x,a)$.
\end{itemize}
\end{assumption}

Since $\sX$ is locally compact separable metric space, there exists a nested sequence of compact sets $\{K_n\}$ such that $K_n \subset \intr K_{n+1}$ and $\sX = \bigcup_{n=1}^{\infty} K_n$ \citet[Lemma 2.76, p. 58]{AlBo06}.

\begin{lemma} \label{nlemma0}
For any compact subset $K$ of $\sX$ and for any $\varepsilon>0$, there exists a compact subset $K_{\varepsilon}$ of $\sX$ such that
\begin{align}
\sup_{(x,a) \in K\times\sA} \int_{K_{\varepsilon}^c} w(y) p(dy|x,a) < \varepsilon, \nonumber
\end{align}
where $D^c$ denotes the complement of the set $D$.
\end{lemma}

\proof{Proof.}
We prove the lemma by contradiction. Assume the claim is wrong. Since every compact subset $K$ of $\sX$ is a subset of $K_n$ for some $n$, the negation of the above lemma is equivalent to the following statement: there exists a compact set $K\subset \sX$ and $\varepsilon>0$ such that for all $n\geq1$ we have
\begin{align}
\sup_{(x,a) \in K\times\sA} \int_{K_{n}^c} w(y) p(dy|x,a) \geq \varepsilon. \nonumber
\end{align}
Note that $w$ is integrable with respect to the probability measures in the set $\bigl\{p(\,\cdot\,|x,a): (x,a) \in K\times\sA\bigr\}$ since
\begin{align}
\sup_{(x,a) \in K\times\sA} \int_{\sX} w(y) p(dy|x,a) \leq \alpha \sup_{x\in K} w(x) < \infty. \nonumber
\end{align}
For each $n$, we prove that $\int_{(\intr K_n)^c} w(y) p(dy|x,a)$ is an upper semi-continuous function on $K\times\sA$. Recall that $\int_{\sX} w(y) p(dy|x,a)$ is a continuous function of $(x,a)$. Let $(x_k,a_k) \rightarrow (x,a)$ in $K\times\sA$. Then $p(\,\cdot\,|x_k,a_k) \rightarrow p(\,\cdot\,|x,a)$ weakly and $\int_{\sX} w(y) p(dy|x_k,a_k) \rightarrow \int_{\sX} w(y) p(dy|x,a)$ by our assumption. If we take $f_k=g_k=f=g=w$ in \citet[Theorem 3.3]{Ser82}, this result implies that $\nu_k(\,\cdot\,)  \rightarrow \nu(\,\cdot\,)$ weakly, where
\begin{align}
\nu_k(D) &= \int_{D} w(y) p(dy|x_k,a_k) \nonumber \\
\nu(D) &= \int_{D} w(y) p(dy|x,a), \nonumber
\end{align}
for all $D\in\B(\sX)$. Then, by \citet[Theorem A]{Bar61} we have
\begin{align}
\int_{(\intr K_{n})^c} w(y) p(dy|x,a) &\coloneqq \nu\bigl((\intr K_{n})^c\bigr) \nonumber \\
&\geq \limsup_{k\rightarrow\infty} \nu_k\bigl((\intr K_{n})^c\bigr) \coloneqq \limsup_{k\rightarrow\infty} \int_{(\intr K_{n})^c} w(y) p(dy|x_k,a_k). \nonumber
\end{align}
Hence, $\int_{(\intr K_n)^c} w(y) p(dy|x,a)$ is upper semi-continuous. Since $K\times\sA$ is compact, there exists $(x_n,a_n) \in K\times\sA$ such that
\begin{align}
\sup_{(x,a) \in K\times\sA} \int_{(\intr K_n)^c} w(y) p(dy|x,a) = \int_{(\intr K_n)^c} w(y) p(dy|x_n,a_n). \nonumber
\end{align}
The sequence $\{(x_n,a_n)\}$ (being a sequence in a compact set $K\times\sA$) has an converging subsequence $\{(x_{n_k},a_{n_k})\}$ with the limit $(x,a) \in K\times\sA$. Then, for all $m\geq2$, we have
\begin{align}
\int_{K_{m-1}^c} w(y) p(dy|x,a) &\geq \int_{(\intr K_m)^c} w(y) p(dy|x,a) \nonumber \\
&\geq \limsup_{k\rightarrow\infty} \int_{(\intr K_m)^c} w(y) p(dy|x_{n_k},a_{n_k}) \nonumber \\
&\geq \limsup_{k\rightarrow\infty} \int_{(\intr K_{n_k})^c} w(y) p(dy|x_{n_k},a_{n_k}) \geq \varepsilon, \nonumber
\end{align}
where the third inequality follows from the fact that $(\intr K_m)^c \supset (\intr K_{n_k})^c$ for $k$ sufficiently large.
But this is a contradiction because $w$ is $p(\,\cdot\,|x,a)$ integrable. \Halmos
\endproof

Let $\{\nu_n\}$ be a sequence of probability measures such that for each $n\geq1$, $\nu_n \in \P(K_n^c)$ and
\begin{align}
\gamma_n &\coloneqq \int_{K_n^c} w(x) \nu_n(dx) < \infty,  \label{eq33} \\
\gamma &= \sup_{n} \tau_n \coloneqq \sup_{n} \max \biggl\{0, \phantom{i} \sup_{(x,a) \in \sX\times\sA} \int_{K_n^c} \bigl( \gamma_n - w(y) \bigr) \hspace{5pt}  p(dy|x,a) \biggr\} < \infty \label{eq3}.
\end{align}
For example, such probability measures can be constructed by choosing $x_n \in K_n^c$ such that $w(x_n) < \inf_{x\in K_n^c} w(x) + \frac{1}{n}$ and letting $\nu_n(\,\cdot\,) = \delta_{x_n}(\,\cdot\,)$.

Similar to the finite-state MDP construction in Section~\ref{compact}, we define a sequence of compact-state MDPs, denoted as c-MDP$_n$, to approximate the original model. To this end, for each $n$ let $\sX_n = K_n \cup \{\Delta_n\}$, where $\Delta_n \in K_n^c$ is a so-called pseudo-state. We define the transition probability $p_n$ on $\sX_n$ given $\sX_n\times\sA$ and the one-stage cost function $c_n: \sX_n\times\sA \rightarrow [0,\infty)$ by
\begin{align}
p_n(\,\cdot\,|x,a) &= \begin{cases}
p\bigl(\,\cdot\, \cap K_n |x,a\bigr) + p\bigl(K_n^c|x,a\bigr) \delta_{\Delta_n},   &\text{ if } x\in K_n  \\
\int_{K_n^c} \biggl( p\bigl(\,\cdot\, \cap K_n |z,a\bigr) + p\bigl(K_n^c|z,a\bigr) \delta_{\Delta_n} \biggr) \nu_n(dz) ,  &\text{ if } x=\Delta_n,
\end{cases} \nonumber \\
c_n(x,a) &= \begin{cases}
c(x,a),   &\text{ if } x\in K_n  \\
\int_{K_n^c} c(z,a) \nu_n(dz) ,  &\text{ if } x=\Delta_n. \nonumber
\end{cases}
\end{align}
With these definitions, c-MDP$_n$ is defined as a Markov decision process with the components $\bigl( \sX_n, \sA, p_n, c_n \bigr)$. History spaces, policies, and cost functions are defined in a similar way as in the original model. Let $\Pi_n$, $\Phi_n$, and $\mathbb{F}_n$ denote the set of all policies, randomized stationary policies and deterministic stationary policies of c-MDP$_n$, respectively. For each policy $\pi \in \Pi_n$ and initial distribution $\mu \in \P(\sX_n)$, we denote the cost functions for c-MDP$_n$ by $J_n(\pi,\mu)$ and $V_n(\pi,\mu)$.

To obtain the main result of this section, we introduce, for each $n$, another MDP, denoted by $\overline{\text{MDP}}_n$, with the components $\bigl(\sX,\sA, q_n,b_n)$ where
\begin{align}
q_n(\,\cdot\,|x,a) &= \begin{cases}
p(\,\cdot\,|x,a),   &\text{ if } x\in K_n  \\
\int_{K_n^c} p\bigl(\,\cdot\,|z,a) \nu_n(dz) ,  &\text{ if } x \in K_n^c,
\end{cases} \nonumber \\
b_n(x,a) &= \begin{cases}
c(x,a),   &\text{ if } x\in K_n  \\
\int_{K_n^c} c(z,a) \nu_n(dz) ,  &\text{ if } x \in K_n^c. \nonumber
\end{cases}
\end{align}
For each policy $\pi \in \Pi$ and initial distribution $\mu \in \P(\sX)$, we denote the cost functions for $\overline{\text{MDP}}_n$ by $\bar{J}_n(\pi,\mu)$ and $\bar{V}_n(\pi,\mu)$.

\subsection{Discounted Cost.}\label{noncompact:disc}

In this section we consider \textbf{(Q1)} for the discounted cost criterion with a discount factor $\beta \in (0,1)$. Throughout this section, it is assumed that Assumption~\ref{as1} holds. The following result states that c-MDP$_n$ and $\overline{\text{MDP}}_n$ are equivalent for the discounted cost.

\begin{lemma}\label{lemma1}
We have
\begin{align}
\bar{J}_n^*(x) = \begin{cases}
J_n^*(x),   &\text{ if } x\in K_n  \\
J_n^*(\Delta_n) ,  &\text{ if } x \in K_n^c,
\end{cases}\label{eq50}
\end{align}
where $\bar{J}_n^*$ is the discounted value function of $\overline{\text{MDP}}_n$ and $J_n^*$ is the discounted value function of c-MDP$_n$, provided that there exist optimal deterministic stationary policies for $\overline{\text{MDP}}_n$ and c-MDP$_n$. Furthermore, if, for any deterministic stationary policy $f \in \mathbb{F}_n$, we define $\bar{f}(x)=f(x)$ on $K_n$ and $\bar{f}(x)=f(\Delta_n)$ on $K_n^c$, then
\begin{align}
\bar{J}_n(\bar{f},x) = \begin{cases}
J_n(f,x),   &\text{ if } x\in K_n  \\
J_n(f,\Delta_n) ,  &\text{ if } x \in K_n^c.
\end{cases}\label{eq51}
\end{align}
In particular, if the deterministic stationary policy $f_n^* \in \mathbb{F}_n$ is optimal for c-MDP$_n$, then its extension $\bar{f}_n^*$ to $\sX$ is also optimal for $\overline{\text{MDP}}_n$.
\end{lemma}

\proof{Proof.}
The proof of (\ref{eq51}) is a consequence of the following facts: $b_n(x,a) = b_n(y,a)$ and $q_n(\,\cdot\,|x,a) = q_n(\,\cdot\,|y,a)$ for all $x,y \in K_n^c$ and $a\in\sA$. In other words, $K_n^c$ in $\overline{\text{MDP}}_n$ behaves like the pseudo state $\Delta_n$ in c-MDP$_n$ when $\bar{f}$ is applied to $\overline{\text{MDP}}_n$.

Let $\overline{\rF}_n$ denote the set of all deterministic stationary policies in $\rF$ which are obtained by extending policies in $\rF_n$ to $\sX$. If we can prove that $\min_{f \in \rF} \bar{J}_n(f,x) = \min_{f \in \overline{\rF}_n} \bar{J}_n(f,x)$ for all $x \in \sX$, then (\ref{eq50}) follows from (\ref{eq51}). Let $f \in \rF \setminus \overline{\rF}_n$. We have two cases: (i) $\bar{J}_n(f,z) = \bar{J}_n(f,y)$ for all $z,y \in K_n^c$ or (ii) there exists $z,y \in K_n^c$ such that $\bar{J}_n(f,z) < \bar{J}_n(f,y)$.

For the case (i), if we define the deterministic Markov policy $\pi^0$ as $\pi^0 = \{f_0,f,f,\ldots\}$, where $f_0(x) = f(z)$ on $K_n^c$ for some fixed $z \in K_n^c$ and $f_0(x) = f(x)$ on $K_n$, then using the expression
\begin{align}
\bar{J}_n(\pi^0,x) = b_n(x,f_0(x)) + \beta \int_{\sX} \bar{J}_n(f,x') q_n(dx'|x,f_0(x)), \label{eq52}
\end{align}
it is straightforward to show that $\bar{J}_n(\pi^0,x) = \bar{J}_n(f,x)$ on $K_n$ and $\bar{J}_n(\pi^0,x) = \bar{J}_n(f,z)$ on $K_n^c$. Therefore, $\bar{J}_n(\pi^0,x) = \bar{J}_n(f,x)$ for all $x\in \sX$ since $\bar{J}_n(f,x) = \bar{J}_n(f,z)$ for all $x \in K_n^c$. For all $t\geq1$ define the deterministic Markov policy $\pi^t$ as $\pi^t = \{f_0,\pi^{t-1}\}$. Analogously, one can prove that $\bar{J}_n(\pi^t,x) = \bar{J}_n(\pi^{t+1},x)$ for all $x\in \sX$. Since $\bar{J}_n(\pi^t,x) \rightarrow \bar{J}_n(f_0,x)$ as $t\rightarrow\infty$, we have $\bar{J}_n(f_0,x) = \bar{J}_n(f,x)$ for all $x\in\sX$, where $f_0 \in \overline{\rF}_n$.

For the second case, if we again consider the deterministic Markov policy $\pi^0 = \{f_0,f,f,\ldots\}$, then by (\ref{eq52}) we have $\bar{J}_n(\pi^0,y) = \bar{J}_n(f,z) < \bar{J}_n(f,y)$. Since $\min_{f \in \rF} \bar{J}_n(f,y) \leq \bar{J}_n(\pi^0,y)$,  this completes the proof. \Halmos
\endproof

For each $n$, let us define $w_n$ by letting $w_n(x) = w(x)$ on $K_n$ and $w_n(x) = \int_{K_n^c} w(z) \nu_n(dz) \eqqcolon \gamma_n$ on $K_n^c$. Hence, $w_n \in B(\sX)$ by (\ref{eq33}).

\begin{lemma} \label{nlemma1}
For all $n$ and $x\in\sX$, the components of $\overline{\text{MDP}}_n$ satisfy the following:
\begin{align}
\sup_{a\in\sA} |b_n(x,a)| &\leq M w_n(x) \label{eq7} \\
\sup_{a\in\sA} \int_{\sX} w_n(y) q_n(dy|x,a) &\leq \alpha w_n(x) + \gamma,\label{eq8}
\end{align}
where $\gamma$ is the constant in (\ref{eq3}).
\end{lemma}

\proof{Proof.}
It is straightforward to prove (\ref{eq7}) by using the definitions of $b_n$ and $w_n$, and the equation (\ref{eq1}). To prove (\ref{eq8}), we have to consider two cases: $x\in K_n$ and $x\in K_n^c$. For the first case, $q_n(\,\cdot\,|x,a) = p(\,\cdot\,|x,a)$, and therefore, we have
\begin{align}
\sup_{a \in \sA} \int_{\sX} w_n(y) p(dy|x,a) &= \sup_{a\in\sA} \biggl\{ \int_{\sX} w(y) p(dy|x,a) + \int_{K_n^c} \bigl( \gamma_n - w(y) \bigr) \hspace{3pt} p(dy|x,a) \biggl\} \nonumber \\
&\leq \sup_{a\in\sA} \int_{\sX} w(y) p(dy|x,a) + \gamma \text{  } \text{ (by (\ref{eq3}))}  \nonumber \\
&\leq \alpha w(x) + \gamma = \alpha w_n(x) + \gamma \text{  } \text{  (as $w_n=w$ on $K_n$)}. \nonumber
\end{align}
For $x\in K_n^c$, we have
\begin{align}
\sup_{a\in\sA} \int_{\sX} w_n(y) q_n(dy|x,a) &= \sup_{a\in \sA} \int_{K_n^c} \biggl( \int_{\sX} w_n(y) p(dy|z,a) \biggr) \nu_n(dz) \nonumber \\
&\leq \int_{K_n^c} \biggl( \sup_{a\in \sA} \int_{\sX} w_n(y) p(dy|z,a) \biggr) \nu_n(dz) \nonumber \\
&\leq \int_{K_n^c} \bigl(\alpha w(z) + \gamma \bigr) \hspace{3pt} \nu_n(dz) \label{aux1} \\
&= \alpha w_n(x) + \gamma, \nonumber
\end{align}
where (\ref{aux1}) can be proved following the same arguments as for the case $x\in K_n$. This completes the proof. \Halmos
\endproof

Note that if we define $c_{n,0}(x) = 1 + \sup_{a\in\sA} |b_n(x,a)|$ and $c_{n,t}(x) = \sup_{a\in\sA} \int_{\sX} c_{n,t-1}(y) q_n(dy|x,a)$, by (\ref{eq7}) and (\ref{eq8}), and an induction argument, we obtain (see \citet[p. 46]{HeLa99})
\begin{align}
c_{n,t}(x) \leq L w_n(x) \alpha^t + L \gamma \sum_{j=0}^{t-1} \alpha^j \text{   } \text{  for all $x\in \sX$}, \label{aux2}
\end{align}
where $L = 1+M$. Let $\beta_0 > \beta$ be such that $\alpha \beta_0 < 1$ and let $C_n: \sX \rightarrow [1,\infty)$ be defined by
\begin{align}
C_n(x) = \sum_{t=0}^{\infty} \beta_0^t c_{n,t}(x). \nonumber
\end{align}
Then, for all $x\in \sX$, by (\ref{aux2}) we have
\begin{align}
C_n(x) \coloneqq \sum_{t=0}^{\infty} \beta_0^t c_{n,t}(x) &\leq \frac{L}{1-\beta_0 \alpha} w_n(x) + \frac{L \beta_0}{(1-\beta_0)(1-\beta_0 \alpha)} \gamma \nonumber \\
&\coloneqq L_1 w_n(x) + L_2. \label{eq9}
\end{align}
Hence $C_n \in B(\sX)$ as $w_n \in B(\sX)$. Moreover, for all $(x,a) \in \sX\times\sA$, $C_n$ satisfies (see \citet[p. 45]{HeLa99})
\begin{align}
\int_{\sX} C_n(y) q_n(dy|x,a) &= \sum_{t=0}^{\infty} \beta_0^t \int_{\sX} c_{n,t}(y) q_n(dy|x,a) \nonumber \\
&\leq \sum_{t=0}^{\infty} \beta_0^t c_{n,t+1}(x) \nonumber \\
&\leq \frac{1}{\beta_0} \sum_{t=0}^{\infty} \beta_0^t c_{n,t}(x) = \alpha_0 C_n(x), \nonumber  
\end{align}
where $\alpha_0 \coloneqq \frac{1}{\beta_0}$ and $\alpha_0 \beta <1$ since $\beta_0 > \beta$. Therefore, for all $x \in \sX$, components of $\overline{\text{MDP}}_n$ satisfy
\begin{align}
\sup_{a\in\sA} |b_n(x,a)| &\leq C_n(x) \label{eq11} \\
\sup_{a\in\sA} \int_{\sX} C_n(y) q_n(dy|x,a) &\leq \alpha_0 C_n(x) \label{eq12}.
\end{align}

Since $C_n \in B(\sX)$, the Bellman optimality operator $\overline{T}_n$ of $\overline{\text{MDP}}_n$ maps $B^l(\sX)$ into $B^l(\sX)$ and is given by
\begin{align}
\overline{T}_n u(x) &= \inf_{a\in\sA} \biggl[ b_n(x,a) + \beta \int_{\sX} u(y) q_n(dy|x,a) \biggr] \nonumber \\
&= \begin{cases}
\inf_{a\in\sA} \bigl[ c(x,a) + \beta \int_{\sX} u(y) p(dy|x,a) \bigr],   &\text{ if } x\in K_n  \\
\inf_{a\in\sA} \int_{K_n^c} \bigl[ c(z,a) + \beta \int_{\sX} u(y) p(dy|z,a) \bigr] \nu_n(dz),  &\text{ if } x \in K_n^c.
\end{cases} \nonumber 
\end{align}
Then successive approximations to the discounted value function of $\overline{\text{MDP}}_n$ are given by $v_n^0 = 0$ and $v_n^{t+1} = \overline{T}_n v_n^{t}$ ($t\geq1$). Since $\alpha_0 \beta <1$, it can be proved as in \citet[Theorem 8.3.6, p. 47]{HeLa99} and \citet[(8.3.34), p. 52]{HeLa99} that
\begin{align}
|v_n^t(x)| , |\bar{J}_n^*(x)| &\leq  \frac{C_n(x)}{1-\sigma_0}  \text{ } \text{ for all $x$,}\label{eq13} \\
\|v_n^t - \bar{J}_n^*\|_{C_n} &\leq  \frac{\sigma_0^t}{1-\sigma_0^t}, \label{eq14}
\end{align}
where $\sigma_0 = \beta \alpha_0 < 1$.

Similar to $v_n^t$, let us define $v^0 = 0$ and $v^{t+1} = T v^t$, where $T: B_w(\sX) \rightarrow B_w(\sX)$, the Bellman optimality operator for the original MDP, is given by
\begin{align}
T u(x) = \inf_{a\in\sA} \biggl[ c(x,a) + \beta \int_{\sX} u(y) p(dy|x,a) \biggr]. \nonumber
\end{align}
Then, again by \citet[Theorem 8.3.6, p. 47]{HeLa99} and \citet[(8.3.34), p. 52]{HeLa99} we have
\begin{align}
|v^t(x)| , |J^*(x)| &\leq M \frac{w(x)}{1-\sigma} \text{ } \text{ for all $x$,} \label{eq15} \\
\|v^t - J^*\|_{w} &\leq M \frac{\sigma^t}{1-\sigma}, \label{eq16}
\end{align}
where $\sigma = \beta \alpha < 1$.

\begin{lemma}\label{lemma2}
For any compact set $K\subset \sX$, we have
\begin{align}
\lim_{n\rightarrow\infty} \sup_{x\in K} |v_n^t(x) - v^t(x)| = 0 \label{eq17}
\end{align}
for all $t\geq1$.
\end{lemma}

\proof{Proof.}
We prove (\ref{eq17}) by induction on $t$. For $t=1$, the claim trivially holds since any compact set $K\subset \sX$ is inside $K_n$ for sufficiently large $n$, and therefore, $b_n = c$ on $K$ for sufficiently large $n$ (recall $v_n^0 = v^0 = 0$). Assume the claim is true for $t\geq1$. Fix any compact set $K$. Recall the definition of compact subsets $K_{\varepsilon}$ of $\sX$ in Lemma~\ref{nlemma0}. By definition of $q_n$, $b_n$, and $w_n$, there exists $n_0\geq1$ such that for all $n\geq n_0$, $q_n = p$, $b_n = c$, and $w_n = w$ on $K$. With these observations, for each $n\geq n_0$ we have
\begin{align}
\sup_{x\in K} |v_n^{t+1}(x) - v^{t+1}(x)|
&= \sup_{x \in K} \biggl| \inf_{\sA} \biggl[ c(x,a) + \beta \hspace{-5pt} \int_{\sX} v_n^t(y) p(dy|x,a) \biggr] - \min_{\sA} \biggl[ c(x,a) + \beta \hspace{-5pt} \int_{\sX} v^t(y) p(dy|x,a) \biggr] \biggr|  \nonumber \\
&\leq \beta \sup_{(x,a) \in K\times\sA} \biggl| \int_{\sX} v_n^t(y) p(dy|x,a) - \int_{\sX} v^t(y) p(dy|x,a) \biggr| \nonumber \\
&= \beta \sup_{(x,a) \in K\times\sA} \biggl| \int_{K_{\varepsilon}} \bigl(v_n^t(y) - v^t(y)\bigr) \hspace{3pt} p(dy|x,a) + \int_{K_{\varepsilon}^c} \bigl(v_n^t(y) - v^t(y)\bigr) \hspace{3pt} p(dy|x,a)\biggr| \nonumber\\
&\leq \beta \biggl\{ \sup_{x \in K_{\varepsilon}} |v_n^t(x) - v^t(x)| + \sup_{(x,a) \in K\times\sA} \biggl| \int_{K_{\varepsilon}^c} \bigl(v_n^t(y) - v^t(y)\bigr) \hspace{3pt} p(dy|x,a)\biggr| \biggr\} \nonumber
\end{align}
Note that we have $|v^t| \leq M \frac{w}{1-\sigma}$ by (\ref{eq15}). Since $w_n \leq \gamma_{\max} w$, where $\gamma_{\max} \coloneqq \max\{1,\gamma\}$, we also have
$|v_n^t| \leq \frac{L_1 \gamma_{\max} w + L_2}{1-\sigma_0} \leq \frac{(L_1 \gamma_{\max}  + L_2) w}{1-\sigma_0}$ by (\ref{eq9}) and (\ref{eq13}) (as $w\geq1$).
Let us define
\begin{align}
R \coloneqq \frac{L_1 \gamma_{\max} + L_2}{1-\sigma_0}+\frac{M}{1-\sigma}. \nonumber
\end{align}
Then by Lemma~\ref{nlemma0} we have
\begin{align}
\sup_{x\in K} |v_n^{t+1}(x) - v^{t+1}(x)| &\leq \beta \sup_{x \in K_{\varepsilon}} |v_n^t(x) - v^t(x)| + \beta  R \varepsilon. \nonumber
\end{align}
Since the first term converges to zero as $n\rightarrow\infty$ by the induction hypothesis, and $\varepsilon$ is arbitrary, the claim is true for $t+1$. This completes the proof. \Halmos
\endproof

The following theorem states that the discounted value function of $\overline{\text{MDP}}_n$ converges to the discounted value function of the original MDP uniformly on each compact set $K\subset\sX$.

\begin{theorem}\label{mainthm1}
For any compact set $K \subset \sX$ we have
\begin{align}
\lim_{n\rightarrow\infty} \sup_{x\in K} | \bar{J}_n^*(x) - J^*(x) | = 0. \label{eq20}
\end{align}
\end{theorem}

\proof{Proof.}
Fix any compact set $K\subset\sX$. Since $w$ is continuous and therefore bounded on $K$, it is sufficient to prove $\lim_{n\rightarrow\infty} \sup_{x\in K} \frac{| \bar{J}_n^*(x) - J^*(x) |}{w(x)}$. Let $n$ be chosen such that $K \subset K_n$, and so, $w_n=w$ on $K$. Then we have
\begin{align}
\sup_{x\in K} \frac{| \bar{J}_n^*(x) - J^*(x) |}{w(x)}
&\leq \sup_{x\in K} \frac{| \bar{J}_n^*(x) - v_n^t(x) |}{w(x)} + \sup_{x\in K} \frac{| v_n^t(x) - v^t(x) |}{w(x)} + \sup_{x\in K} \frac{| v^t(x) - J^*(x) |}{w(x)} \nonumber \\
&\leq \sup_{x\in K} \frac{| \bar{J}_n^*(x) - v_n^t(x) |}{C_n(x)} \frac{C_n(x)}{w(x)} + \sup_{x\in K} \frac{| v_n^t(x) - v^t(x) |}{w(x)} +
M \frac{\sigma^t}{1-\sigma^t} \text{  } \text{ (by (\ref{eq16}))} \nonumber \\
&\leq \sup_{x\in K} \frac{| \bar{J}_n^*(x) - v_n^t(x) |}{C_n(x)} \frac{(L_1 w_n(x)+L_2)}{w(x)} + \sup_{x\in K} \frac{| v_n^t(x) - v^t(x) |}{w(x)} + \frac{M \sigma^t}{1-\sigma^t} \text{} \text{ (by (\ref{eq9}))} \nonumber \\
&\leq (L_1 + L_2) \sup_{x\in K} \frac{| \bar{J}_n^*(x) - v_n^t(x) |}{C_n(x)} + \sup_{x\in K} \frac{| v_n^t(x) - v^t(x) |}{w(x)} +
\frac{M \sigma^t}{1-\sigma^t} \text{} \text{ ($w_n=w$ on $K$)} \nonumber \\
&\leq (L_1+L_2) \frac{\sigma_0^t}{1-\sigma_0} + \sup_{x\in K} \frac{| v_n^t(x) - v^t(x) |}{w(x)} +
\frac{M \sigma^t}{1-\sigma^t} \text{  } \text{ (by (\ref{eq14}))}. \nonumber
\end{align}
Since $w\geq1$ on $\sX$, $\sup_{x\in K} \frac{| v_n^t(x) - v^t(x) |}{w(x)} \rightarrow 0$ as $n\rightarrow\infty$ for all $t$ by Lemma~\ref{lemma2}. Hence, the last expression can be made arbitrarily small. This completes the proof. \Halmos
\endproof

In the remainder of this section, we use the above results and Theorem~\ref{compact:mainthm1} to compute a near optimal policy for the original MDP. It is straightforward to check that for each $n$, c-MDP$_n$ satisfies the assumptions in Theorem~\ref{compact:mainthm1}. Let $\{\varepsilon_n\}$ be a sequence of positive real numbers such that $\lim_{n\rightarrow\infty} \varepsilon_n = 0$.

By Theorem~\ref{compact:mainthm1}, for each $n\geq1$, there exists a deterministic stationary policy $f_{n} \in \mathbb{F}_n$, obtained from the finite state approximations of c-MDP$_n$, such that
\begin{align}
\sup_{x \in \sX_n}| J_n(f_{n},x) - J_n^*(x) | \leq \varepsilon_n, \nonumber 
\end{align}
where for each $n$, finite-state models are constructed replacing $\bigl( \sZ,\sA,p,c \bigr)$ with the components $\bigl( \sX_n,\sA,p_n,c_n \bigr)$ of c-MDP$_n$ in Section~\ref{compact}. By Lemma~\ref{lemma1}, for each $n\geq1$ we also have
\begin{align}
\sup_{x \in \sX}| \bar{J}_n(f_{n},x) - \bar{J}_n^*(x) | \leq \varepsilon_n, \label{eq22}
\end{align}
where, with an abuse of notation, we also denote the extended (to $\sX$) policy by $f_n$. Let us define operators $\bar{R}_{n}: B_{C_n}(\sX) \rightarrow B_{C_n}(\sX)$ and $R_n: B_w(\sX) \rightarrow B_w(\sX)$ by
\begin{align}
\bar{R}_n u(x)  &= \begin{cases}
c(x,f_{n}(x)) + \beta \int_{\sX} u(y) p(dy|x,f_{n}(x)),   &\text{ if } x\in K_n  \\
\int_{K_n^c} \bigl[ c(z,f_{n}(z)) + \beta \int_{\sX} u(y) p(dy|z,f_{n}(z)) \bigr] \nu_n(dz),  &\text{ if } x \in K_n^c,
\end{cases} \nonumber \\
R_n u(x) &= c(x,f_{n}(x)) + \beta \int_{\sX} u(y) p(dy|x,f_{n}(x)). \nonumber 
\end{align}
By \citet[Remark 8.3.10, p. 54]{HeLa99}, $\bar{R}_n$ is a contraction operator with modulus $\sigma_0$ and $R_n$ is a contraction operator with modulus $\sigma$. Furthermore, the fixed point of $\bar{R}_n$ is $\bar{J}_n(f_{n},x)$ and the fixed point of $R_n$ is $J(f_{n},x)$. For each $n\geq1$, let us define $\bar{u}_n^0 = u_n^0 = 0$ and $\bar{u}_n^{t+1} = \bar{R}_n \bar{u}_n^{t}$, $u_n^{t+1} = R_n u_n^{t}$ ($t\geq1$). One can prove that (see the proof of \citet[Theorem 8.3.6, p. 51]{HeLa99})
\begin{align}
|\bar{u}_n^t(x)| , |\bar{J}_n(f_{n},x)| &\leq \frac{C_n(x)}{1-\sigma_0} \nonumber \\
\| \bar{u}_n^t - \bar{J}_n(f_{n},\,\cdot\,)\|_{C_n} &\leq  \frac{\sigma_0^t}{1-\sigma_0} \nonumber \\
|u_n^t(x)| , |J(f_{n},x)| &\leq M \frac{w(x)}{1-\sigma} \nonumber \\
\| u_n^t - J(f_{n},\,\cdot\,)\|_{w} &\leq M \frac{\sigma^t}{1-\sigma}. \nonumber 
\end{align}

\begin{lemma}\label{lemma3}
For any compact set $K \subset \sX$, we have
\begin{align}
\lim_{n\rightarrow\infty} \sup_{x \in K} |\bar{u}_n^t(x) - u_n^t(x)| = 0. \nonumber 
\end{align}
\end{lemma}

\proof{Proof.}
The lemma can be proved using the same arguments as in the proof of Lemma~\ref{lemma2} and so we omit the details. \Halmos
\endproof

\begin{lemma}\label{lemma4}
For any compact set $K \subset \sX$, we have
\begin{align}
\lim_{n\rightarrow\infty} \sup_{x \in K} |\bar{J}_n(f_{n},x) - J(f_{n},x)| = 0. \label{eq30}
\end{align}
Indeed, this is true for all sequences of policies in $\rF$.
\end{lemma}

\proof{Proof.}
The lemma can be proved using the same arguments as in the proof of Theorem~\ref{mainthm1}. \Halmos
\endproof

The following theorem is the main result of this section which states that the true cost functions of the policies obtained from finite state models converge to the value function of the original MDP. Hence, to obtain a near optimal policy for the original MDP, it is sufficient to compute the optimal policy for the finite state model that has sufficiently large number of grid points.

\begin{theorem}\label{mainthm2}
For any compact set $K \subset \sX$, we have
\begin{align}
\lim_{n\rightarrow\infty} \sup_{x\in K} |J(f_{n},x) - J^*(x)| &= 0. \nonumber 
\intertext{Therefore,}
\lim_{n\rightarrow\infty} |J(f_{n},x) - J^*(x)| &= 0 \text{  } \text{ for all $x \in \sX$}. \nonumber 
\end{align}
\end{theorem}

\proof{Proof.}
The result follows from (\ref{eq20}), (\ref{eq22}), and (\ref{eq30}). \Halmos
\endproof

\subsection{Average Cost.}\label{sec3}

In this section we obtain approximation results, analogous to Theorems~\ref{mainthm1} and \ref{mainthm2}, for the average cost criterion. To do this, we impose some new assumptions on the components of the original MDP in addition to Assumption~\ref{as1}.
These assumptions are the unbounded counterpart of Assumption~\ref{compact:as2}. With the exception of Assumption~\ref{as2}-(j), versions of these assumptions were imposed in \citet{Veg03}, \citet{GoHe95}, and \citet{JaNo06} to study the existence of the solution to the Average Cost Optimality Equality (ACOE) and Inequality (ACOI). In what follows, for any finite signed measure $\vartheta$ and measurable function $h$ on $\sX$, we let $\vartheta(h) \coloneqq \int_{\sX} h(x) \vartheta(dx)$ and
\begin{align}
\|\vartheta\|_w \coloneqq \sup_{\|g\|_w\leq 1} \biggl| \int_{\sX} g(x) \vartheta(dx) \biggr|. \nonumber
\end{align}
Here $\|\vartheta\|_w$ is called the $w$-norm of $\vartheta$.

\begin{assumption}\label{as2}
Suppose Assumption~\ref{as1} holds with item (b) and (\ref{eq2}) replaced by conditions (j) and (e) below, respectively. In addition, there exist a probability measure $\eta$ on $\sX$ and a positive measurable function $\phi:\sX\times\sA \rightarrow (0,\infty)$ such that for all $(x,a) \in \sX\times\sA$
\begin{itemize}
\item [(e)] $\int_{\sX} w(y) p(dy|x,a) \leq \alpha w(x) + \eta(w) \phi(x,a)$, where $\alpha \in (0,1)$.
\item [(f)] $p(D|x,a) \geq \eta(D) \phi(x,a)$ for all $D \in \B(\sX)$.
\item [(g)] The weight function $w$ is $\eta$-integrable, i.e., $\eta(w) < \infty$.
\item [(h)] For each $n\geq1$, $\inf_{(x,a)\in K_n\times\sA}\phi(x,a)>0$.
\item [(j)] The stochastic kernel $p(\,\cdot\,|x,a)$ is continuous in $(x,a)$ with respect to the $w$-norm.
\end{itemize}
\end{assumption}

Throughout this section, it is assumed that Assumption~\ref{as2} holds. Conditions (e), (f), and (g) of Assumption~\ref{as2} are unbounded counterparts of conditions (d) and (e) in Assumption~\ref{compact:as2}. Recall that condition (e) corresponds to the so-called `drift inequality' and condition (f) corresponds to the so-called `minorization' condition which guarantee the geometric ergodicity of Markov chains induced by stationary policies (see \citet{HeLa99,MeTw93} and references therein). These assumptions are quite general for studying average cost problems with unbounded one-stage costs. In addition, they are proper for the approximation problem in the sense that it can be shown that if the original problem satisfies these, then the reduced models constructed in the sequel satisfy similar conditions. 
There is only one minor difference between Assumption~\ref{as2}-(f) and the standard minorization condition: in the literature $\phi$ is in general required to be nonnegative instead of positive.

Note that although Assumption~\ref{as2}-(j) seems to be restrictive, it is weaker than the assumptions imposed in the literature for studying approximation of average cost problems with unbounded cost (see \citet{DuPr14}). Indeed, it is assumed in \citet{DuPr14} that the transition probability $p$ is Lipschitz continuous in $(x,a)$ with respect to $w$-norm. The reason for imposing such a strong condition on the transition probability is to obtain convergence rate for the approximation problem. Since we do not aim to provide rate of convergence result in this section, it is natural to impose continuity instead of Lipschitz continuity of the transition probability. However, it does not seem possible to replace continuity with respect to the $w$-norm by a weaker convergence notion. One reason is that with a weaker continuity notion it is not possible to prove that the transition probability of c-MDP$_n$ is continuous with respect to the total variation distance, which is needed if one wants to use Theorem~\ref{compact:thm6} and cannot be relaxed as explained in Remark~\ref{remarktotal}.

Analogous with Theorem~\ref{compact:thm4}, the following theorem is a consequence of \citet[Theorems 3.3]{Veg03}, \citet[Lemma 3.4]{GoHe95} (see also \citet[Proposition 10.2.5, p. 126]{HeLa99}), and \citet[Theorem 3]{JaNo06}, which also holds with Assumption~\ref{as2}-(j) replaced by Assumption~\ref{as1}-(b).

\begin{theorem}\label{thm1}
For each $f \in \rF$, the stochastic kernel $p(\,\cdot\,|x,f(x))$ is positive Harris recurrent with unique invariant probability measure $\mu_f$. Furthermore, $w$ is $\mu_f$-integrable, and therefore, $\rho_f  \coloneqq \int_{\sX} c(x,f) \mu_f(dx) <~\infty$.
There exist positive real numbers $R$ and $\kappa < 1$ such that
\begin{align}
\sup_{f\in\rF} \| p^t(\,\cdot\,|x,f(x)) - \mu_f \|_{w} \leq R w(x) \kappa^t \label{eq40}
\end{align}
for all $x\in\sX$, where $R$ and $\kappa$ continuously depend on $\alpha$, $\eta(w)$, and $\inf_{f\in\rF} \eta(\phi(y,f(y)))$.
Finally, there exist $f^{*} \in \rF$ and  $h^{*} \in B_w(\sX)$ such that the triplet $(h^{*},f^{*},\rho_{f^{*}})$ satisfies the average cost optimality equality (ACOE), and therefore,
\begin{align}
\inf_{\pi \in \Pi} V(\pi,x)  \coloneqq  V^{*}(x) = \rho_{f^{*}}, \nonumber
\end{align}
for all $x\in\sX$.
\end{theorem}
Note that (\ref{eq40}) implies that for each $f \in \rF$, the average cost is given by $V(f,x) = \int_{\sX} c(y,f(y)) \mu_f(dy)$ for all $x\in\sX$ (instead of $\mu_f$-a.e.); that is, the average cost is independent of the initial point.

Recall that $V_n$ and $\bar{V}_n$ denote the average costs of c-MDP$_n$ and $\overline{\text{MDP}}_n$, respectively. The value functions for average cost are denoted analogously to the discounted cost case. Similar to Lemma~\ref{lemma1}, the following result states that MDP$_n$ and $\overline{\text{MDP}}_n$ are not too different for the average cost.

\begin{lemma}\label{lemma5}
Suppose Theorem~\ref{thm1} holds for $\overline{\text{MDP}}_n$ and Theorem~\ref{compact:thm4} holds for MDP$_n$. Then we have
\begin{align}
\bar{V}_n^*(x) = \begin{cases}
V_n^*(x),   &\text{ if } x\in K_n  \\
V_n^*(\Delta_n) ,  &\text{ if } x \in K_n^c.
\end{cases}\label{eq53}
\end{align}
Furthermore, if, for any deterministic stationary policy $f \in \mathbb{F}_n$, we define $\bar{f}(x)=f(x)$ on $K_n$ and $\bar{f}(x)=f(\Delta_n)$ on $K_n^c$, then
\begin{align}
\bar{V}_n(\bar{f},x) = \begin{cases}
V_n(f,x),   &\text{ if } x\in K_n  \\
V_n(f,\Delta_n) ,  &\text{ if } x \in K_n^c.
\end{cases}\label{eq54}
\end{align}
In particular, if the deterministic stationary policy $f_n^* \in \mathbb{F}_n$ is optimal for MDP$_n$, then its extension $\bar{f}_n^*$ to $\sX$ is also optimal for $\overline{\text{MDP}}_n$.
\end{lemma}

\proof{Proof.}
Let the triplet $(h_n^{*},f_n^{*},\rho^n_{f_n^{*}})$ satisfy the ACOE for c-MDP$_n$, so that $f_n^{*}$ is an optimal policy and $\rho^n_{f_n^{*}}$ is the average value function for c-MDP$_n$. It is straightforward to show that the triplet $(\bar{h}_n^{*},\bar{f}_n^{*},\rho^n_{f_n^{*}})$ satisfies the ACOE for $\widetilde{\text{MDP}}_n$, where
\begin{align}
\bar{h}_n^*(x) &= \begin{cases}
h_n^*(x),   &\text{ if } x\in K_n  \\
h_n^*(\Delta_n),  &\text{ if } x\in K_n^c,
\end{cases}\nonumber \\
\intertext{and}
\bar{f}_n^*(x) &= \begin{cases}
f_n^*(x),   &\text{ if } x\in K_n  \\
f_n^*(\Delta_n),  &\text{ if } x\in K_n^c.
\end{cases}\nonumber
\end{align}
By \citet[Lemma 5.2]{GoHe95} (see also \citet[Section 5.2]{HeLa96}), this implies that $\bar{f}_n^{*}$ is an optimal stationary policy for $\overline{\text{MDP}}_n$ with cost function $\rho^n_{f_n^{*}}$. This completes the proof of the first part.

For the second part, let $f\in\rF_n$ with an unique invariant probability measure $\mu_{f} \in \P(\sX_n)$ and let $\bar{f} \in \rF$ denote its extension to $\sX$ with an unique invariant probability measure $\mu_{\bar{f}}$. It can be proved that
\begin{align}
\mu_f(\,\cdot\,) = \mu_{\bar{f}}(\,\cdot\,\cap K_n) + \mu_{\bar{f}}(K_n^c) \delta_{\Delta_n}(\,\cdot\,). \nonumber
\end{align}
Then we have
\begin{align}
\bar{V}_n(f,x) &= \int_{\sX} b_n(x,\bar{f}(x)) \mu_{\bar{f}}(dx) \nonumber  \\
&= \int_{K_n} c_n(x,\bar{f}(x)) \mu_{\bar{f}}(dx) + \mu_{\bar{f}}(K_n^c) c_n(\Delta_n,\bar{f}(\Delta_n)) \nonumber \\
&= \int_{\sX_n} c_n(x,f(x)) \mu_f(dx) \nonumber \\
&= V_n(f,x). \nonumber
\end{align}
This completes the proof. \Halmos
\endproof

By Lemma~\ref{lemma5}, in the remainder of this section we need only consider $\overline{\text{MDP}}_n$ in place of MDP$_n$. Later we will show that Theorem~\ref{thm1} holds for $\overline{\text{MDP}}_n$ for $n$ sufficiently large and that Theorem~\ref{compact:thm4} holds for c-MDP$_n$ for all $n$.

Recall the definition of constants $\gamma_n$ and $\tau_n$ from (\ref{eq33}) and (\ref{eq3}). For each $n\geq1$, we define $\phi_n: \sX\times\sA \rightarrow (0,\infty)$ and $\varsigma_n \in \R$ as
\begin{align}
\phi_n(x,a) &\coloneqq \begin{cases}
\phi(x,a), &\text{ if } x \in K_n \\
\int_{K_n^c} \phi(y,a) \nu_n(dy), &\text{ if } x \in K_n^c,
\end{cases} \nonumber \\
\varsigma_n &\coloneqq \int_{K_n^c} w(y) \eta(dy). \nonumber
\end{align}
Since $\eta(w) < \infty$ and $\tau_n$ can be made arbitrarily small by properly choosing $\nu_n$, we assume, without loss of generality, the following.

\begin{assumption}\label{as-aux}
The sequence of probability measures $\{\nu_n\}$ is chosen such that the following holds
\begin{align}
\lim_{n\rightarrow\infty} (\tau_n + \varsigma_n) = 0. \label{eq37}
\end{align}
\end{assumption}

Let $\alpha_n \coloneqq \alpha + \varsigma_n + \tau_n$.

\begin{lemma}
For all $n$ and $(x,a) \in \sX\times\sA$, the components of $\overline{\text{MDP}}_n$ satisfy the following:
\begin{align}
\sup_{a\in\sA} |b_n(x,a)| &\leq M w_n(x) \nonumber \\*
\int_{\sX} w_n(y) q_n(dy|x,a) &\leq \alpha_n w_n(x) + \eta(w_n) \phi_n(x,a),\label{neq10} \\
q_n(D|x,a) &\geq \eta(D) \phi_n(x,a) \text{  } \text{ for all $D\in\B(\sX)$}.\nonumber
\end{align}
\end{lemma}

\proof{Proof.}
The proof of the first inequality follows from Assumption~\ref{as2} and definitions of $b_n$ and $w_n$. To prove the remaining two inequalities, we have to consider the cases $x\in K_n$ and $x\in K_n^c$ separately.

Let $x\in K_n$, and therefore, $q_n(\,\cdot\,|x,a) = p(\,\cdot\,|x,a)$. The second inequality holds since
\begin{align}
\int_{\sX} w_n(y) p(dy|x,a) &= \int_{\sX} w(y) p(dy|x,a) + \int_{K_n^c} \bigl(\gamma_n - w(y) \bigr) \hspace{3pt} p(dy|x,a) \nonumber \\
&\leq \int_{\sX} w(y) p(dy|x,a) + \tau_n \nonumber \\
&\leq \alpha w(x) + \eta(w) \phi(x,a) + \tau_n \nonumber \\
&\leq \alpha w_n(x) + \eta(w_n) \phi_n(x,a) + \varsigma_n \phi_n(x,a) + \tau_n \text{  } \text{  (as $w_n=w$ and $\phi_n=\phi$ on $K_n$)} \nonumber \\
&\leq \alpha_n w_n(x) + \eta(w_n) \phi_n(x,a), \text{   } \text{ (as $\phi_n\leq1$ and $w_n\geq1$)}.\nonumber
\end{align}
For the last inequality, for all $D\in \B(\sX)$, we have
\begin{align}
q_n(D|x,a) &= p(D|x,a) \geq \eta(D) \phi(x,a) = \eta(D) \phi_n(x,a) \text{   } \text{   (as $\phi_n=\phi$ on $K_n$)}.\nonumber
\end{align}
Hence, inequalities hold for $x\in K_n$.

For $x\in K_n^c$, we have
\begin{align}
\int_{\sX} w_n(y) q_n(dy|x,a) &= \int_{K_n^c} \biggl( \int_{\sX} w_n(y) p(dy|z,a) \biggr) \nu_n(dz) \nonumber \\
&\leq \int_{K_n^c} \bigl( \alpha w(z) + \eta(w_n) \phi(x,a) + \varsigma_n \phi(x,a) + \tau_n \bigr) \hspace{3pt} \nu_n(dz) \label{aux3} \\
&= \alpha w_n(x) + \eta(w_n) \phi_n(x,a) + \varsigma_n \phi_n(x,a) + \tau_n \nonumber \\
&\leq \alpha_n w_n(x) + \eta(w_n) \phi_n(x,a), \text{   } \text{ (since $\phi_n\leq1$ and $w_n\geq1$)}\nonumber
\end{align}
where (\ref{aux3}) can be obtained following the same arguments as for the case $x\in K_n$. The last inequality holds for $x\in K_n^c$ since
\begin{align}
q_n(D|x,a) &= \int_{K_n^c} p(D|z,a) \nu_n(dz) \nonumber \\
&\geq \int_{K_n^c} \eta(D) \phi(z,a) \nu_n(dz) \nonumber \\
&= \eta(D) \phi_n(x,a). \nonumber
\end{align}
This completes the proof. \Halmos
\endproof

We note that by (\ref{eq37}), there exists $n_0\geq1$ such that $\alpha_n < 1$ for $n\geq n_0$. Hence, for each $n\geq n_0$, Theorem~\ref{thm1} holds for $\overline{\text{MDP}}_n$ with $w$ replaced by $w_n$ for some $R_n > 0$ and $\kappa_n \in (0,1)$, and we have $R_{\max} \coloneqq \sup_{n\geq n_0} R_n < \infty$ and $\kappa_{\max} \coloneqq \sup_{n \geq n_0} \kappa_n <~ 1$.

In the remainder of this section, it is assumed that $n\geq n_0$.

\begin{lemma}\label{lemma6}
Let $g:\sX\times\sA\rightarrow \R$ be any measurable function such that $\sup_{a\in\sA} |g(x,a)| \leq M_g w(x)$ for some $M_g \in \R$. Then, for all $t\geq1$ and any compact set $K\subset\sX$ we have
\begin{align}
\sup_{(y,f) \in K\times\rF} \biggl| \int_{\sX} g_n(x,f(x)) q_n^t(dx|y,f(y)) - \int_{\sX} g(x,f(x)) p^t(dx|y,f(y)) \biggr| \rightarrow 0 \nonumber
\end{align}
as $n\rightarrow\infty$, where $g_n(x,a) = g(x,a)$ on $K_n\times\sA$ and $g_n(x,a) = \int_{K_n^c} g(z,a) \nu_n(dz)$ on $K_n^c\times\sA$.
\end{lemma}

\proof{Proof.}
We will prove the lemma by induction. Fix any compact set $K\subset \sX$.  We note that in the inequalities below, we repeatedly use the fact $\phi , \phi_n \leq1$ without explicitly referring to this fact. Recall the definition of the compact subsets $K_{\varepsilon}$ of $\sX$ in Lemma~\ref{nlemma0} and the constant $\gamma_{\max} = \max\{1,\gamma\}$. Note that $\sup_{a\in \sA} |g_n(x,a)| \leq M_g w_n(x) \leq M_g \gamma_{\max} w(x)$ for all $x\in\sX$.

The claim holds for $t=1$ by the following argument:
\begin{align}
&\sup_{(y,f) \in K\times\rF} \biggl| \int_{\sX} g_n(x,f(x)) q_n(dx|y,f(y)) - \int_{\sX} g(x,f(x)) p(dx|y,f(y)) \biggr| \nonumber \\
&\phantom{xxx}=\sup_{(y,f) \in K\times\rF} \biggl| \int_{\sX} g_n(x,f(x)) p(dx|y,f(y)) - \int_{\sX} g(x,f(x)) p(dx|y,f(y)) \biggr| \text{ } \text{ (for $n$ sufficiently large)} \nonumber \\
&\phantom{xxx}=\sup_{(y,f) \in K\times\rF} \biggl| \int_{K_{\varepsilon}^c} g_n(x,f(x)) p(dx|y,f(y)) - \int_{K_{\varepsilon}^c} g(x,f(x)) p(dx|y,f(y)) \biggr| \text{ } \text{ (for $n$ sufficiently large)} \nonumber \\
&\phantom{xxx}\leq M_g(1+\gamma_{\max}) \varepsilon \nonumber,
\end{align}
where the last inequality follows from Lemma~\ref{nlemma0}. Since $\varepsilon$ is arbitrary, the result follows.

Assume the claim is true for $t\geq1$. Let us define
$l_f(z) := \int_{\sX} g(x,f(x)) p^t(dx|z,f(z))$ and $l^n_f(z) := \int_{\sX} g_n(x,f(x)) q_n^t(dx|z,f(z))$. By recursively applying the inequalities in Assumption~\ref{as2}-(e) and in (\ref{neq10}) we obtain
\begin{align}
\sup_{f\in \rF} |l_f(z)| &\leq M_g \alpha^t w(z) + M_g \eta(w) \sum_{j=0}^{t-1} \alpha^j \nonumber \\
\intertext{and}
\sup_{f\in \rF} |l^n_f(z)| &\leq M_g \alpha_n^t w_n(z) + M_g \eta(w_n) \sum_{j=0}^{t-1} \alpha_n^j \nonumber \\
&\leq M_g \alpha_{\max}^t \gamma_{\max} w(z) + M_g \eta(w) \gamma_{\max} \sum_{j=0}^{t-1} \alpha_{\max}^j, \nonumber
\end{align}
where $\alpha_{\max} \coloneqq \sup_{n\geq n_0} \alpha_n < 1$. Then we have
\begin{align}
&\sup_{(y,f) \in K\times\rF} \biggl| \int_{\sX} g_n(x,f(x)) q_n^{t+1}(dx|y,f(y)) - \int_{\sX} g(x,f(x)) p^{t+1}(dx|y,f(y)) \biggr| \nonumber \\*
&\phantom{xxxxxxxx}=\sup_{(y,f) \in K\times\rF}  \biggl| \int_{\sX} l^n_f(z) q_n(dz|y,f(y)) - \int_{\sX} l_f(z) p(dz|y,f(y)) \biggr| \nonumber \\
&\phantom{xxxxxxxx}= \sup_{(y,f) \in K\times\rF} \biggl| \int_{\sX} l^n_f(z) p(dz|y,f(y)) - \int_{\sX} l_f(z) p(dz|y,f(y)) \biggr| \text{ (for $n$ sufficiently large)} \nonumber \\
&\phantom{xxxxxxxx}\leq \sup_{(y,f) \in K\times\rF} \biggl| \int_{K_{\varepsilon}^c} l^n_f(z) p(dz|y,f(y)) - \int_{K_{\varepsilon}^c} l_f(z) p(dz|y,f(y)) \biggr| + \hspace{-10pt} \sup_{(z,f) \in K_{\varepsilon}\times\rF} |l^n_f(z) - l_f(z)| \nonumber \\
&\phantom{xxxxxxxx}\leq R \varepsilon + \hspace{-5pt} \sup_{(z,f) \in K_{\varepsilon}\times\rF} |l^n_f(z) - l_f(z)|  \label{eq60},
\end{align}
where $R$ is given by
\begin{align}
R \coloneqq M_g \biggl( \alpha^t + \alpha_{\max}^t \gamma_{\max} + \eta(w)\sum_{j=0}^{t-1} \alpha^j + \eta(w) \gamma_{\max} \sum_{j=0}^{t-1} \alpha_{\max}^j \biggl) \nonumber
 \end{align}
and the last inequality follows from Lemma~\ref{nlemma0}. Since the claim holds for $t$ and $K_{\varepsilon}$, the second term in (\ref{eq60}) goes to zero as $n\rightarrow\infty$. Since $\varepsilon$ is arbitrary, the result follows. \Halmos
\endproof

In the remainder of this section the above results are used to compute a near optimal policy for the original MDP. Let $\{\varepsilon_n\}$ be a sequence of positive real numbers converging to zero.

For each $f \in \rF$, let $\mu^n_f$ denote the unique invariant probability measure of the transition kernel $q_n(\,\cdot\,|x,f(x))$ and let $\rho^n_f$ denote the associated average cost; that is, $\rho^n_f \coloneqq \bar{V}_n(f,x) = \int_{\sX} b_n(y,f(y)) \mu^n_f(dy)$ for all initial points $x\in\sX$. Therefore, the value function of $\overline{\text{MDP}}_n$, denoted by $\bar{V}^*_n$, is given by $V^*_n(x) = \inf_{f\in\rF} \rho_f^n$, i.e., it is constant on $\sX$.

Before making the connection with Theorem~\ref{compact:thm6}, we prove the following result.

\begin{lemma}\label{lemma8}
The transition probability $p_n$ of c-MDP$_n$ is continuous in $(x,a)$ with respect to the total variation distance.
\end{lemma}

\proof{Proof.}
To ease the notation, we define $M(\sX_n)$, $M(\sX)$, and $M_w(\sX)$ as the subsets of $B(\sX_n)$, $B(\sX)$, and $B_w(\sX)$, respectively, whose elements have (corresponding) norm less than one. Let $(x_k,a_k) \rightarrow (x,a)$ in $\sX_n\times\sA$. Since the pseudo state $\Delta_n$ is isolated and $K_n$ is compact, we have two cases: (i) $x_k = x = \Delta_n$ for all $k$ large enough, or (ii) $x_k \rightarrow x$ in $K_n$.

For the first case we have
\begin{align}
\|p_n(\,\cdot\,|\Delta_n,a_k) - p_n(\,\cdot\,|\Delta_n,a)\|_{TV}
&= \sup_{g \in M(\sX_n)} \biggl| \int_{\sX_n} g(y) p_n(dy|\Delta_n,a_k) - \int_{\sX_n} g(y) p_n(dy|\Delta_n,a) \biggr| \nonumber \\
&\leq \sup_{g \in M(\sX)} \biggl| \int_{\sX} g(y) q_n(dy|\Delta_n,a_k) - \int_{\sX} g(y) q_n(dy|\Delta_n,a) \biggr| \label{eqq1} \\
&= \sup_{g \in M(\sX)} \biggl| \int_{K_n^c} \biggl( \int_{\sX} g(y) p(dy|z,a_k) - \int_{\sX} g(y) p(dy|z,a) \biggr) \nu_n(dz) \biggr| \nonumber \\
&\leq \int_{K_n^c} \sup_{g\in M(\sX)} \biggl| \int_{\sX} g(y) p(dy|z,a_k) - \int_{\sX} g(y) p(dy|z,a) \biggr| \nu_n(dz) \nonumber \\
&\leq \int_{K_n^c} \sup_{g\in M_w(\sX)} \biggl| \int_{\sX} g(y) p(dy|z,a_k) - \int_{\sX} g(y) p(dy|z,a) \biggr| \nu_n(dz) \nonumber \\
&= \int_{K_n^c} \|p(\,\cdot\,|z,a_k) - p(\,\cdot\,|z,a)\|_w \nu_n(dz) \label{eqq2},
\end{align}
where (\ref{eqq1}) follows since if for any $g \in M(\sX_n)$ we define $\bar{g} = g$ on $K_n$ and $\bar{g} = g(\Delta_n)$ on $K_n^c$, then we have $\bar{g} \in M(\sX)$ and $\int_{\sX_n} g(y) p_n(dy|x,a) = \int_{\sX} \bar{g}(y) q_n(dy|x,a)$ for all $(x,a) \in \sX_n\times\sA$. Note that we have
\begin{align}
\sup_{g\in M_w(\sX)} \biggl| \int_{\sX} g(y) p(dy|z,a_k) - \int_{\sX} g(y) p(dy|z,a) \biggr|
&\leq \int_{\sX} w(y) p(dy|z,a_k) + \int_{\sX} w(y) p(dy|z,a)  \nonumber \\
&\leq 2 \bigl( \alpha + \eta(w) \bigr) w(z)  \nonumber
\end{align}
by Assumption~\ref{as2}-(e), $\phi\leq1$, and $w\geq1$. Since $w$ (restricted to $K_n^c$) is $\nu_n$-integrable, by the dominated convergence theorem (\ref{eqq2}) goes to zero as $k\rightarrow\infty$.

For the second case we have
\begin{align}
\|p_n(\,\cdot\,|x_k,a_k) - p_n(\,\cdot\,|x,a)\|_{TV} &= \hspace{-5pt} \sup_{g \in M(\sX_n)} \biggl| \int_{\sX_n} g(y) p_n(dy|x_k,a_k) - \int_{\sX_n} g(y) p_n(dy|x,a) \biggr| \nonumber \\
&\leq \sup_{g \in M(\sX)} \biggl| \int_{\sX} g(y) q_n(dy|x_k,a_k) - \int_{\sX} g(y) q_n(dy|x,a) \biggr| \nonumber \\
&= \sup_{g \in M(\sX)} \biggl| \int_{\sX} g(y) p(dy|x_k,a_k) - \int_{\sX} g(y) p(dy|x,a) \biggr| \text{ } \text{ (since $x_k,x \in K_n$)}\nonumber \\
&\leq \sup_{g \in M_w(\sX)} \biggl| \int_{\sX} g(y) p(dy|x_k,a_k) - \int_{\sX} g(y) p(dy|x,a) \biggr| \nonumber \\
&= \|p(\,\cdot\,|x_k,a_k) - p(\,\cdot\,|x,a)\|_w. \nonumber
\end{align}
By Assumption~\ref{as2}-(j) the last term goes to zero as $k\rightarrow\infty$. \Halmos
\endproof

Thus we obtain that for each $n\geq 1$, c-MDP$_n$ satisfies the assumption in Theorem~\ref{compact:thm6} for
\begin{align}
\zeta(\,\cdot\,) &= \eta(\,\cdot\, \cap K_n) + \eta(K_n^c) \delta_{\Delta_n}(\,\cdot\,), \nonumber \\
\theta(x,a) &= \begin{cases}
\phi(x,a),   &\text{ if } x\in K_n  \\
\int_{K_n^c} \phi(y,a) \nu_n(dy) ,  &\text{ if } x = \Delta_n,
\end{cases} \nonumber
\end{align}
and some $\lambda \in (0,1)$, where the existence of $\lambda$ follows from Assumption~\ref{as2}-(h) and the fact that $\phi > 0$.

Consequently, there exists a deterministic stationary policy $f_n \in \rF_n$, obtained from the finite state approximations of c-MDP$_n$, such that
\begin{align}
\sup_{x \in \sX_n} |V_n(f_n,x) - V_n^*(x)| \leq \varepsilon_n, \label{eq34}
\end{align}
where finite-state models are constructed replacing $\bigl( \sZ,\sA,p,c \bigr)$ with the components $\bigl( \sX_n,\sA,p_n,c_n \bigr)$ of c-MDP$_n$ in Section~\ref{compact}. By Lemma~\ref{lemma5}, we also have
\begin{align}
|\rho^n_{f_n} - \bar{V}^*_n| \leq \varepsilon_n, \label{eq35}
\end{align}
where, by an abuse of notation, we also denote the policy extended to $\sX$ by $f_n$.

\begin{lemma}\label{lemma7}
We have
\begin{align}
\sup_{f\in\rF} |\rho^n_f - \rho_f| \rightarrow 0 \label{eq36}
\end{align}
as $n\rightarrow\infty$.
\end{lemma}

\proof{Proof.}
Fix any compact set $K\subset \sX$. For any $t\geq1$ and $y\in K$, we have
\begin{align}
\sup_{f\in\rF}  |\rho^n_f - \rho_f| &= \sup_{f\in\rF} \biggl| \int_{\sX} b_n(x,f(x)) \mu^n_f(dx) - \int_{\sX} c(x,f(x)) \mu_f(dx) \biggr| \nonumber \\*
&\leq \sup_{f \in \rF} \biggl | \int_{\sX} b_n(x,f(x)) \mu_f^n(dx) - \int_{\sX} b_n(x,f(x)) q_n^t(dx|y,f(y)) \biggr | \nonumber \\
&\phantom{xxxxxxxxxxx}+ \sup_{f \in \rF} \biggl | \int_{\sX} b_n(x,f(x)) q_n^t(dx|y,f(y)) - \int_{\sX} c(x,f(x)) p^t(dx|y,f(y)) \biggr | \nonumber \\
&\phantom{xxxxxxxxxxxxxxxxx}+ \sup_{f \in \rF} \biggl | \int_{\sX} c(x,f(x)) p^t(dx|y,f(y)) - \int_{\sX} c(x,f(x)) \mu_f(dx) \biggr | \nonumber \\
&\leq M R_{\max} w(y) \kappa_{\max}^t + M R w(y) \kappa^t +\nonumber \\
&\phantom{xxxxxxxxxxx} \sup_{(y,f)\in K\times\rF} \biggl | \int_{\sX} b_n(x,f(x)) q_n^t(dx|y,f(y)) - \int_{\sX} c(x,f(x)) p^t(dx|y,f(y)) \biggr | , \nonumber
\end{align}
where the last inequality follows from Theorem~\ref{thm1}-(ii) and (\ref{eq1}) in Assumption~\ref{as1}. The result follows from Lemma~\ref{lemma6}. \Halmos
\endproof

\begin{theorem}\label{mainthm3}
The value function of $\overline{\text{MDP}}_n$ converges to the value function of the original MDP, i.e., $|\bar{V}^*_n - V^*| \rightarrow 0$, as $n\rightarrow\infty$.
\end{theorem}

\proof{Proof.}
Since
\begin{align}
|\bar{V}^*_n - V^*| &= |\inf_{f\in\rF} \rho^n_f - \inf_{f\in\rF} \rho_f| \leq \sup_{f\in\rF} |\rho^n_f - \rho_f|, \nonumber
\end{align}
the result follows from Lemma~\ref{lemma7}. \Halmos
\endproof

The following is the main result of this section which states that the true average cost of the policies $f_n$ obtained from finite state approximations of c-MDP$_n$ converges to the average value function $V^*$ of the original MDP.

\begin{theorem}\label{mainthm4}
We have $|\rho_{f_n} - V^*| \rightarrow 0$, as $n\rightarrow\infty$.
\end{theorem}

\proof{Proof.}
We have
\begin{align}
|\rho_{f_n} - V^*| &\leq |\rho_{f_n} - \rho^n_{f_n}| + |\rho^n_{f_n} - \bar{V}^*_n| + |\bar{V}^*_n - V^*| \nonumber \\
&\leq \sup_{f\in\rF} |\rho_f - \rho^n_f| + \varepsilon_n + |\bar{V}^*_n - V^*| \text{      } \text{ (by (\ref{eq35}))} \nonumber
\end{align}
The result follows from Lemma~\ref{lemma7} and Theorem~\ref{mainthm3}. \Halmos
\endproof

\section{Discretization of the Action Space.}\label{act dist}

For computing near optimal policies using well known algorithms, such as value iteration, policy iteration, and $Q$-learning, the action space must be finite. In this section, we show that, as a pre-processing step, the action space can taken to be finite if it has sufficiently large number of points for accurate approximation. Throughout this section, it is assumed that Assumption~\ref{as1} holds for the discounted cost and Assumption~\ref{as2} holds for the average cost.

It was shown in \citet{SaLiYu13-2} and \citet{SaYuLi16} that any MDP with (infinite) compact action space can be well approximated by an MDP with finite action space under assumptions that are satisfied by c-MDP$_n$, for both the discounted cost and the average cost cases. Specifically, let $d_{\sA}$ denote the metric on $\sA$. Since $\sA$ is compact, one can find a sequence of finite subsets $\{\Lambda_k\}$ of $\sA$ such that for all $k$
\begin{align}
\min_{\hat{a}\in\Lambda_k} d_{\sA}(a,\hat{a}) < 1/k \text{  } \text{for all $a\in\sA$}. \nonumber
\end{align}
We define c-MDP$_{n,k}$ as the Markov decision process having the components $\bigl\{ \sX_n,\Lambda_k,p_n,c_n \bigr\}$ and we let $\rF_n(\Lambda_k)$ denote the set of all deterministic stationary policies for c-MDP$_{n,k}$. Note that $\rF_n(\Lambda_k)$ is the set of policies in $\rF_n$ taking values only in $\Lambda_k$. Therefore, in a sense, c-MDP$_{n,k}$ and c-MDP$_n$ can be viewed as the same MDP, where the former has constraints on the set of policies.
For each $n$ and $k$, by an abuse of notation, let $f_n^*$ and $f_{n,k}^*$ denote the optimal stationary policies of c-MDP$_n$ and c-MDP$_{n,k}$, respectively, for both the discounted and average costs. Then \citet[Theorem 3.2]{SaYuLi16} and \citet[Theorem 3.2]{SaLiYu13-2} show that for all $n$, we have
\begin{align}
\lim_{k\rightarrow\infty} J_{n}(f^*_{n,k},x) &= J_n(f_n^*,x) \coloneqq J_n^*(x) \nonumber \\
\lim_{k\rightarrow\infty} V_{n}(f^*_{n,k},x) &= V_n(f_n^*,x), \coloneqq V_n^*(x) \nonumber
\end{align}
for all $x\in\sX_n$. In other words, the discounted and average value functions of c-MDP$_{n,k}$ converge to the discounted and average value functions of c-MDP$_n$ as $k\rightarrow\infty$. We note that although \citet[Theorem 3.2]{SaYuLi16} and \citet[Theorem 3.2]{SaLiYu13-2} are proved for nonnegative one-stage cost function, it is straightforward to check that these theorems are also valid for any real valued one-stage cost function.

\begin{theorem}
\label{discrete-act}
For any $x \in \sX$, there exists a subsequence $\{k_n\}$ such that
\begin{align}
\lim_{n\rightarrow\infty} J(f^*_{n,k_n},x) &= J^*(x) \nonumber \\
\lim_{n\rightarrow\infty} V(f^*_{n,k_n},x) &= V^*(x), \nonumber
\end{align}
where $f^*_{n,k_n} \in \rF(\Lambda_{k_n})$ is the optimal stationary policy of c-MDP$_{n,k_n}$.
\end{theorem}

\proof{Proof.}
Let us fix $x\in\sX$. For $n$ sufficiently large (so $x \in K_n$), we choose $k_n$ such that $|J_{n}(f^*_{n,k_n},x)-J_n(f_n^*,x)|< 1/n$ (or $|V_{n}(f^*_{n,k_n},x)-V_n(f_n^*,x)|< 1/n$ for the average cost). We note that if $\sA$ is a compact subset of a finite dimensional Euclidean space, then by using \citet[Theorems 4.1 and 4.2]{SaLiYu13-2} one can obtain an explicit expression for $k_n$ in terms of $n$ under further continuity conditions on $c$ and $p$. By Lemmas~\ref{lemma4} and \ref{lemma7}, we have $|\bar{J}_n(f^*_{n,k_n},x) - J(f^*_{n,k_n},x)|\rightarrow 0$ and $|\bar{V}_n(f^*_{n,k_n},x) - V(f^*_{n,k_n},x)|\rightarrow0$ as $n\rightarrow\infty$, where again by an abuse of notation, the policies extended to $\sX$ are also denoted by $f^*_{n,k_n}$.
Since $\bar{J}_n(f^*_{n,k_n},x) = J_n(f^*_{n,k_n},x)$ and $\bar{V}_n(f^*_{n,k_n},x) = V_n(f^*_{n,k_n},x)$, using Theorems~\ref{mainthm1} and \ref{mainthm3} one can immediately obtain
\begin{align}
\lim_{n\rightarrow\infty} J(f^*_{n,k_n},x) &= J^*(x) \nonumber \\
\lim_{n\rightarrow\infty} V(f^*_{n,k_n},x) &= V^*(x). \nonumber
\end{align}
\Halmos
\endproof

Theorem~\ref{discrete-act} implies that before discretizing the state space to compute the near optimal policies, one can discretize, without loss of generality, the action space $\sA$ in advance on a finite grid using sufficiently large number of grid points.

\section{Rate of Convergence Analysis for Compact-State MDPs.}
\label{compact:rateconv}

In this section we consider \textbf{(Q2)} for MDPs with compact state space; that is, we derive an upper bound on the performance loss due to discretization in terms of the cardinality of the set $\sZ_n$ (i.e., number of grid points) . To do this, we will impose some new assumptions on the components of the MDP in addition to Assumptions~\ref{compact:as1} and \ref{compact:as2}. First, we present some definitions that are needed in the development.

For each $g \in C_b(\sZ)$, let
\begin{align}
\|g\|_{\Lip} \coloneqq \sup_{(z,y)\in\sZ\times\sZ} \frac{|g(z)-g(y)|}{d_{\sZ}(z,y)}. \nonumber
\end{align}
If $\|g\|_{\Lip}$ is finite, then $g$ is called Lipschitz continuous with Lipschitz constant $\|g\|_{\Lip}$. $\Lip(\sZ)$ denotes the set of all Lipschitz continuous functions on $\sZ$, i.e.,
\begin{align}
\Lip(\sZ) \coloneqq \{g \in C_b(\sZ): \|g\|_{\Lip} < \infty \} \nonumber
\end{align}
and $\Lip(\sZ,K)$ denotes the set of all $g \in \Lip(\sZ)$ with $\|g\|_{\Lip} \leq K$. The \emph{Wasserstein distance of order $1$} \citet[p. 95]{Vil09} between two probability measures $\zeta$ and $\xi$ over $\sZ$ is defined as
\begin{align}
W_1(\zeta,\xi) \coloneqq \sup \biggl\{\biggl|\int_{\sZ} g d\zeta - \int_{\sZ} g d\xi\biggr|: g \in \Lip(\sZ,1)\biggr\}. \nonumber
\end{align}
$W_1$ is also called the \emph{Kantorovich-Rubinstein distance}. It is known that if $\sZ$ is compact, then $W_1(\zeta,\xi) \leq \diam(\sZ) \|\zeta - \xi\|_{TV}$; see \citet[Theorem 6.15, p. 103]{Vil09}. For compact $\sZ$, the Wasserstein distance of order $1$ is weaker than total variation distance. Furthermore, for compact $\sZ$, the Wasserstein distance of order $1$ metrizes the weak topology on the set of probability measures $\P(\sZ)$ (see \citet[Corollary 6.13, p. 97]{Vil09}) which also implies that convergence in this sense is weaker than setwise convergence.

In this section we impose the following supplementary assumptions in addition to Assumption~\ref{compact:as1} and Assumption~\ref{compact:as2}.

\begin{assumption}
\label{compact:as3}
\begin{itemize}
\item [  ]
\item[(g)] The one-stage cost function $c$ satisfies $c(\,\cdot\,,a) \in \Lip(\sZ,K_1)$ for all $a \in \sA$ for some $K_1$.
\item[(h)] The stochastic kernel $p$ satisfies $W_1\bigl(p(\,\cdot\,|z,a),p(\,\cdot\,|y,a)\bigr) \leq K_2 d_{\sZ}(z,y)$ for all $a \in \sA$ for some $K_2$.
\item[(j)] $\sZ$ is an infinite compact subset of $\R^d$ for some $d\geq1$, equipped with the Euclidean norm.
\end{itemize}
\end{assumption}

We note that Assumption~\ref{compact:as3}-(j) implies the existence of a constant $\alpha>0$ and finite
subsets $\sZ_n\subset\sZ$ with cardinality $n$ such that
\begin{align}
\max_{z\in\sZ}\min_{y\in \sZ_n} d_{\sZ}(z,y)\leq \alpha (1/n)^{1/d}  \label{compact:quantcof}
\end{align}
for all $n$, where $d_{\sZ}$ is the Euclidean distance on $\sZ$.
In the remainder of this section, we replace $\sZ_n$ defined in Section~\ref{compact}
with $\sZ_n$ satisfying (\ref{compact:quantcof}) in order to derive \emph{explicit} bounds on the approximation error in terms of the cardinality of $\sZ_n$.

\subsection{Discounted Cost.}
\label{compact:discrate}

Assumptions~\ref{compact:as1} and \ref{compact:as3} are imposed throughout this section. Additionally, we assume that $K_2 \beta < 1$. The last assumption is the key to prove the next result which states that the value function $J^{*}$ of the original MDP for the discounted cost is in $\Lip(\sZ)$. Although this result is known in the literature (see \citet{Hin05}), we give a short proof for the sake of completeness using a simple application of the value iteration algorithm.

\begin{theorem}\label{compact:lipcont}
The value function $J^{*}$ for the discounted cost is in $\Lip(\sZ,K)$, where $K= K_1 \frac{1}{1-\beta K_2}$.
\end{theorem}

\proof{Proof.}
Let $u \in \Lip(\sZ,K)$ for some $K>1$. Then $g = \frac{u}{K} \in \Lip(\sZ,1)$ and therefore, for all $a \in \sA$ and $z,y \in \sZ$ we have
\begin{align}
\biggl | \int_{\sZ} u(x) p(dx|z,a) - \int_{\sZ} u(x) p(dx|y,a) \biggr | &= K \biggl | \int_{\sZ} g(x) p(dx|z,a) - \int_{\sZ} g(x) p(dx|y,a) \biggr | \nonumber \\
&\leq K W_1\bigl(p(\,\cdot\,|z,a), p(\,\cdot\,|y,a)\bigr) \leq K K_2 d_{\sZ}(z,y), \nonumber
\end{align}
by Assumption~\ref{compact:as3}-(h). Hence, the contraction operator $T$ defined in (\ref{aux6}) maps $u \in \Lip(\sZ,K)$ to $Tu \in \Lip(\sZ,K_1+\beta K K_2)$, since, for all $z,y \in \sZ$
\begin{align}
| Tu(z) - Tu(y) | &\leq \max_{a\in\sA} \biggl \{ |c(z,a) - c(y,a)| + \beta \biggl | \int_{\sZ} u(x) p(dx|z,a) - \int_{\sZ} u(x) p(dx|y,a) \biggr | \biggr \}\nonumber \\
&\leq K_1 d_{\sZ}(z,y) + \beta K K_2 d_{\sZ}(z,y) = \bigl(K_1 + \beta K K_2\bigr) d_{\sZ}(z,y). \nonumber
\end{align}
Now we apply $T$ recursively to obtain the sequence $\{T^n u\}$ by letting $T^n u = T (T^{n-1} u )$, which converges to the value function $J^{*}$ by the Banach fixed point theorem. Clearly, by induction we have for all $n\geq1$
\begin{align}
T^n u \in \Lip(\sZ,K_n), \nonumber
\end{align}
where $K_n = K_1 \sum_{i=0}^{n-1} (\beta K_2)^i + K (\beta K_2)^n$. If we choose $K < K_1$, then $K_n \leq K_{n+1}$ for all $n$ and therefore, $K_n \uparrow K_1 \frac{1}{1-\beta K_2}$ since $K_2 \beta < 1$. Hence, $T^n u \in \Lip(\sZ,K_1 \frac{1}{1-\beta K_2})$ for all $n$, and therefore, $J^{*} \in \Lip(\sZ,K_1 \frac{1}{1-\beta K_2})$ since $\Lip(\sZ,K_1 \frac{1}{1-\beta K_2})$ is closed with respect to the sup-norm $\|\,\cdot\,\|$. \Halmos
\endproof

The following theorem is the main result of this section. Recall that the policy $\hf_n \in \rF$ is obtained by extending the optimal policy $f_n^{*}$ of MDP$_n$ to $\sZ$.

\begin{theorem}\label{compact:mainthm3}
We have
\begin{align}
\|J(\hf_n,\,\cdot\,) - J^{*}\| \leq \frac{\tau(\beta,K_2)K_1 \frac{1}{1-\beta K_2}+\frac{2K_1}{1-\beta}}{1-\beta} 2 \alpha (1/n)^{1/d}, \nonumber
\end{align}
where $\tau(\beta,K_2) = (2+\beta)\beta K_2 + \frac{\beta^2+4\beta+2}{(1-\beta)^2}$ and $\alpha$ is the coefficient in (\ref{compact:quantcof}).
\end{theorem}

\proof{Proof.}
To prove the theorem, we obtain upper bounds on the expressions derived in Section~\ref{compact:sec2sub1} in terms of the cardinality $n$ of $\sZ_n$. The proof of Theorem~\ref{compact:mainthm1} gives
\begin{align}
\hspace{-2pt}\| J(\hf_n,\,\cdot\,) - J^* \| \leq \frac{\| T_{\hf_n} J^* - \hat{T}_{\hf_n} J^* \| + (1+\beta) \| \hJ^*_n - J^* \|}{1-\beta}. \nonumber
\end{align}
To prove the theorem we upper bound $\| T_{\hf_n} J^* - \hat{T}_{\hf_n} J^* \|$ and $\| \hJ^*_n - J^* \|$ in terms $n$.
For the first term we have
\begin{align}
&\| T_{\hf_n} J^* - \hat{T}_{\hf_n} J^* \| = \sup_{z\in\sZ} \bigl | T_{\hf_n} J^*(z) - \hat{T}_{\hf_n} J^*(z) | \nonumber \\
\small
&\phantom{xxx}\leq \sup_{z\in\sZ} \int \biggl| \scalebox{0.96}{$c(z,\hf_n(z))$} + \beta \int_{\sZ} \scalebox{0.96}{$J^{*}(y) p(dy|z,\hf_n(z))$} - \scalebox{0.96}{$c(x,\hf_n(x))$} - \beta \int_{\sZ} \scalebox{0.96}{$J^{*}(y) p(dy|x,\hf_n(x))$} \biggr| \scalebox{0.96}{$\nu_{n,i_n(z)}(dx)$}\nonumber \\
\normalsize
&\phantom{xxx}\leq\sup_{z\in\sZ} \int \biggl(K_1 d_{\sZ}(x,z) + \beta \biggl| \int_{\sZ} \hspace{-5pt} J^{*}(y) p(dy|z,\hf_n(z)) - \hspace{-5pt} \int_{\sZ} \hspace{-5pt} J^{*}(y) p(dy|x,\hf_n(z)) \biggr| \biggr) \nu_{n,i_n(z)}(dx) \nonumber \\
&\phantom{xxxxxxxxxxxxxxxxxxxxxxxxxxxxx}\text{ (since $\hf_n(x)=\hf_n(z)$ for all  $x\in\S_{n,i_n(z)}$)} \nonumber \\
&\phantom{xxx}\leq \sup_{z\in\sZ} \int (K_1 + \beta \|J^{*}\|_{\Lip} K_2) d_{\sZ}(x,z) \nu_{n,i_n(z)}(dx) \nonumber \\
&\phantom{xxx}\leq (K_1 + \beta \|J^{*}\|_{\Lip} K_2) \max_{i \in \{1,\ldots,n\}} \diam(\S_{n,i}) \nonumber \\
&\phantom{xxx}\leq (K_1 + \beta \|J^{*}\|_{\Lip} K_2) 2\alpha (1/n)^{1/d}. \label{compact:bound1}
\end{align}
For the second term, the proof of Theorem~\ref{compact:thm3} gives
\begin{align}
\| \hJ_n^* - J^* \| \leq \frac{\| \hat{T}_n J^* - F_n J^* \| + (1+\beta) \| J^* - u_n^{*} \|}{1-\beta}. \nonumber
\end{align}
First consider $\| \hat{T}_n J^* - F_n J^* \|$. Define
\begin{align}
l(z,a) &\coloneqq c(z,a) + \beta \int_{\sX} J^{*}(y) p(dy|z,a), \nonumber
\intertext{so that}
J^{*}(z) &= \min_{a\in\sA} l(z,a). \nonumber
\end{align}
It is straightforward to show that $l(\,\cdot\,,a) \in \Lip(\sZ,K_l)$ for all $a \in \sA$, where $K_l = K_1 + \beta \|J^{*}\|_{\Lip} K_2$. By adapting the proof of Lemma~\ref{compact:nlemma3} to the value function $J^{*}$, we obtain
\begin{align}
\| \hat{T}_n J^* -  F_n J^* \| &= \sup_{z\in\sZ} \phantom{x} \biggl | \min_{a\in\sA} \int l(x,a) \nu_{n,i_n(z)}(dx) - \int \min_{a\in\sA} l(x,a) \nu_{n,i_n(z)}(dx) \biggr | \nonumber \\
&\leq \sup_{z \in \sZ} \int \sup_{y \in \S_{n,i_n(z)}} \bigl|l(y,a_i) - J^{*}(y) \bigr| \nu_{n,i_n(z)}(dy) \nonumber \\
&\leq \max_{i \in \{1,\ldots,n\}} \int \sup_{y \in \S_{n,i}} \bigl\{ |l(y,a_i) - l(z_i,a_i)| + |J^{*}(z_i) - J^{*}(y)| \bigr\} \nu_{n,i}(dy) \nonumber \\
&\leq \max_{i \in \{1,\ldots,n\}} \int \sup_{y \in \S_{n,i_n}} \bigl\{ K_l d_{\sZ}(y,z_i) + \|J^{*}\|_{\Lip} d_{\sZ}(z_i,y) \bigr\} \nu_{n,i}(dy) \nonumber \\
&\leq (K_l +  \|J^{*}\|_{\Lip}) \max_{i \in \{1,\ldots,n\}} \diam(\S_{n,i}) \nonumber  \\
&\leq (K_l +  \|J^{*}\|_{\Lip}) 2 \alpha (1/n)^{1/d}. \label{compact:bound2}
\end{align}
For the expression $\|J^{*} - u_n^{*}\|$, by Lemma~\ref{compact:nlemma2} we have
\begin{align}
\| u_n^{*} - J^{*} \| \leq \frac{2}{1-\beta} \inf_{r \in \sZ^{k_n}} \| J^{*} - \Phi_{r} \|, \nonumber
\end{align}
where $\Phi_{r}(z) = \Sigma_{i=1}^{k_n} r_i 1_{S_{n,i}}(z)$, $r = (r_1,\ldots,r_{k_n})$. Since $\|J^{*}\|_{\Lip} < \infty$, we have $\inf_{r \in \sZ^{k_n}} \| J^{*} - \Phi_{r} \| \leq \|J^{*}\|_{\Lip} \max_{i \in \{1,\ldots,n\}} \diam(\S_{n,i}) \leq \|J^{*}\|_{\Lip} 2 \alpha (1/n)^{1/d}$. Hence
\begin{align}
\| u_n^{*} - J^{*} \| \leq \frac{2}{1-\beta} \|J^{*}\|_{\Lip} 2 \alpha (1/n)^{1/d}. \label{compact:bound3}
\end{align}
Hence, by (\ref{compact:bound2}) and (\ref{compact:bound3}) we obtain
\begin{align}
\| \hJ_n^* - J^* \| \leq  \biggl ( \bigl (\beta K_2 + \frac{\beta+3}{(1-\beta)^2} \bigr) \|J^{*}\|_{\Lip} + \frac{K_1}{1-\beta} \biggr) 2 \alpha (1/n)^{1/d}. \label{compact:bound4}
\end{align}
Then, the result follows from (\ref{compact:bound1}) and (\ref{compact:bound4}), and the fact $\|J^{*}\|_{\Lip} \leq K_1 \frac{1}{1-\beta K_2}$. \Halmos
\endproof

\begin{remark}
It is important to point out that if we replace Assumption~\ref{compact:as3}-(h) with the uniform Lipschitz continuity of $p(\,\cdot\,|z,a)$ in $z$ with respect to total variation distance, then Theorem~\ref{compact:mainthm3} remains valid (with possibly different constants in front of the term $(1/n)^{1/d}$). However, in this case, we do not need the assumption $K_2 \beta < 1$.
\end{remark}

\begin{remark}
For the average cost case, instead of assuming from the outset the uniform Lipschitz continuity of $c$ and $p$ in the $z$ variable, we first derive a rate of convergence result in terms of the moduli of continuity of the functions $\omega_c$ and $\omega_p$ in the $z$ variable of $c(z,a)$ and $p(\,\cdot\,|z,a)$, where the total variation distance is used to define $\omega_p$. Then, we state that explicit rate of convergence result can be given if we impose some structural assumptions on $\omega_c$ and $\omega_p$ such as linearity, which corresponds to the uniform Lipschitz continuity of $c(z,a)$ and $p(\,\cdot\,|z,a)$ in $z$. However, this is not the right approach for the discounted cost case as the modulus of continuity function $\omega_p$ is calculated using the Wasserstein distance of order $1$. Indeed, to obtain a similar result as in the average cost case, we must relate $\omega_c$ and $\omega_p$ to the modulus of continuity $\omega_{J^*}$ of the value function $J^*$. This can be established if $\omega_c$ and $\omega_p$ are affine functions (i.e., $\omega_c(r) = K_1 r + L_1$ and $\omega_p(r) = K_2 r + L_2$) using the dual formulation of the Wasserstein distance of order $1$ \cite[Theorem 5.10]{Vil09}:
\begin{align}
W_1(\mu,\nu) = \sup_{\substack{(\psi,\varphi) \in C_b(\sZ)\times C_b(\sZ) \\ \psi(x)-\varphi(y) \leq d_{\sZ}(x,y)}} \biggl| \int_{\sZ} \psi(z) \mu(dz) - \int_{\sZ} \varphi(z) \nu(dz) \biggr|. \nonumber
\end{align}
However, in this situation we can explicitly compute the convergence rate only if $L_1=L_2=0$ which is the uniform Lipschitz continuity case.
\end{remark}

\subsection{Average Cost.}
\label{compact:averate}

In this section, we suppose that Assumptions~\ref{compact:as2} and \ref{compact:as3}-(j) hold. We define the modulus of continuity functions in the $z$ variable of $c(z,a)$ and $p(\,\cdot\,|z,a)$ as follows
\begin{align}
\omega_c(r) &\coloneqq  \sup_{a \in \sA} \sup_{z,y \in \sZ:d_{\sZ}(z,y)\leq r} |c(z,a) - c(y,a)| \nonumber \\
\omega_p(r) &\coloneqq  \sup_{a \in \sA} \sup_{z,y \in \sZ:d_{\sZ}(z,y)\leq r} \|p(\,\cdot\,|z,a)-p(\,\cdot\,|y,a)\|_{TV}. \nonumber
\end{align}
Since $c(z,a)$ and $p(\,\cdot\,|z,a)$ are uniformly continuous, we have $\lim_{r\rightarrow0} \omega_c(r)=0$ and $\lim_{r\rightarrow0} \omega_p(r)=0$. Note that when $\omega_c$ and $\omega_p$ are linear, $c(z,a)$ and $p(\,\cdot\,|z,a)$ are uniformly Lipschitz in $z$. In the remainder of this section, we first derive a rate of convergence result in terms of $\omega_c$ and $\omega_p$. Then, we explicitly compute the convergence rate for the Lipschitz  case as a corollary of this result.

To obtain convergence rates for the average cost, we first prove a rate of convergence result for Lemma~\ref{compact:prop3}. To this end, for each $n\geq1$, let $d_n \coloneqq 2 \alpha (1/n)^{1/d}$, where $\alpha$ is the coefficient in (\ref{compact:quantcof}).

\begin{lemma}\label{compact:ratelemma}
For all $t\geq1$, we have
\begin{align}
\sup_{(y,f) \in \sZ \times \rF} \|p^t(\,\cdot\,|y,f(y)) - q_n^t(\,\cdot\,|y,f(y))\|_{TV} \leq t \omega_p(d_n).  \nonumber
\end{align}
\end{lemma}

\proof{Proof.}
Similar to the proof of Lemma~\ref{compact:prop3}, we use induction. For $t=1$, recalling the proof of Lemma~\ref{compact:prop3}, the claim holds by the following argument:
\begin{align}
\sup_{(y,f)\in \sZ \times \rF} \|p(\,\cdot\,|y,f(y)) - q_n(\,\cdot\,|y,f(y))\|_{TV} &\leq \sup_{y \in \sZ} \sup_{(x,a) \in \S_{n,i_n(y)} \times \sA} \hspace{-15pt}\|p(\,\cdot\,|y,a) - p(\,\cdot\,|x,a)\|_{TV} \nonumber \\
&\leq \omega_p(d_n). \nonumber
\end{align}
Now, assume the claim is true for $t\geq1$. Again recalling the proof of Lemma~\ref{compact:prop3}, we have
\begin{align}
&\sup_{(y,f) \in \sZ \times \rF} \|p^{t+1}(\,\cdot\,|y,f(y)) - q_n^{t+1}(\,\cdot\,|y,f(y))\|_{TV} \leq \sup_{(y,f) \in \sZ\times\rF} \bigl\| p^t(\,\cdot\,|y,f(y)) - q_n^t(\,\cdot\,|y,f(y)) \bigr\|_{TV} \nonumber \\
&\phantom{xxxxxxxxxxxxxxxxxxxxxxxxxxxxxxxxxxxxxxxx} + \hspace{-3pt} \sup_{(z,f) \in \sZ\times\rF} \bigl\| p(\,\cdot\,|z,f(z)) - q_n(\,\cdot\,|z,f(z)) \bigr\|_{TV} \nonumber \\
&\phantom{xxxxxxxxxxxxxxxxxxxxxxxxx}\leq t \omega_p(d_n) + \omega_p(d_n)
=(t+1) \omega_p(d_n). \nonumber
\end{align}
This completes the proof. \Halmos
\endproof

The following theorem is the main result of this section. A somewhat similar result was obtained in \citet[Section 3.5]{Her89}, where identical assumptions are imposed on both the original model and the approximating model (see \citet[Assumption 5.1]{Her89}). Moreover, the approximating transition probability and one-stage cost function are assumed to converge to the original transition probability and one-stage cost function with respect to some rate; that is, $\rho(n) := \sup_{(x,a) \in \sX\times\sA} |b_n(x,a) - c(x,a)|$ and $\pi(n) := \sup_{(x,a) \in \sX\times\sA} \|q_n(\,\cdot\,|x,a) - p(\,\cdot\,|x,a)\|_{TV}$ with $\rho(n), \pi(n) \rightarrow 0$ as $n\rightarrow \infty$. Although our result may appear to be a special case of the results in \citet[Section 3.5]{Her89}, there are several differences:
(i) our assumptions are only imposed for the the original model, and (ii) in \citet[Section 3.5]{Her89} the approximating models do not have finite state space while our approximating models are obtained by extending finite state models to the original state space, thereby, allowing for constructive numerical method to calculate near optimal policies.

Recall that the optimal policy $\tilde{f}^{*}_n$ for $\widetilde{\text{MDP}}_n$ is obtained by extending the optimal policy $f_n^{*}$ for MDP$_n$ to $\sZ$, and $R$ and $\kappa$ are the constants in Theorem~\ref{compact:thm4}.

\begin{theorem}\label{compact:mainthm4}
For all $t\geq1$, we have
\begin{align}
|\rho_{\tilde{f}_n^{*}} - \rho_{f^{*}}| \leq 4 \|c\| R \kappa^t + 2 \omega_c(d_n) + 2 \|c\| t \omega_p(d_n). \nonumber
\end{align}
\end{theorem}

\proof{Proof.}
The proof of Theorem~\ref{compact:thm6} gives
\begin{align}
|\rho_{\tilde{f}_n^{*}} - \rho_{f^{*}}| \leq |\rho_{\tilde{f}_n^{*}} - \hat{\rho}_{\tilde{f}_n^{*}}^n| + |\hat{\rho}_{\tilde{f}_n^{*}}^n - \hat{\rho}_{\hat{f}_n^{*}}^n| + |\hat{\rho}_{\hat{f}_n^{*}}^n - \rho_{f^{*}}|. \nonumber
\end{align}
Hence, to prove the theorem we obtain an upper bounds on the three terms in the sum.
Consider the first term (recall the proof of Lemma~\ref{compact:prop4})
\begin{align}
|\rho_{\tilde{f}_n^{*}} - \hat{\rho}^n_{\tilde{f}_n^{*}}| &\leq \sup_{f \in \rF} |\hat{\rho}_f^n - \rho_f| \nonumber \\
&\leq 2 R \kappa^t \|c\| + \|c\| \sup_{(y,f)\in\sZ\times\rF} \|q_n^t(\,\cdot\,|y,f(y)) - p^t(\,\cdot\,|y,f(y))\|_{TV} \nonumber \\
&\leq 2 R \kappa^t \|c\| + \|c\| t \omega_p(d_n)  \text{ (by Lemma~\ref{compact:ratelemma})}. \label{compact:avebound1}
\end{align}
For the second term, the proof of Lemma~\ref{compact:prop7} gives
\begin{align}
|\hat{\rho}^n_{\tilde{f}_n^{*}} - \hat{\rho}^n_{\hat{f}_n^{*}}| &\leq |\hat{\rho}^n_{\tilde{f}_n^{*}} - \tilde{\rho}^n_{\tilde{f}_n^{*}}| + |\tilde{\rho}^n_{\tilde{f}_n^{*}} - \hat{\rho}^n_{\hat{f}_n^{*}}| \nonumber \\
&\leq \sup_{f\in\rF} |\hat{\rho}^n_{f} - \tilde{\rho}^n_{f}| + |\inf_{f\in\rF}\tilde{\rho}^n_{f} - \inf_{f\in\rF}\hat{\rho}^n_{f}| \nonumber \\
&\leq 2 \sup_{f\in\rF} |\hat{\rho}^n_{f} - \tilde{\rho}^n_{f}| \nonumber \\
&\leq 2 \|b_n - c\| \text{ (see the proof of Lemma~\ref{compact:prop5})} \nonumber \\
&\leq 2 \sup_{(z,a)\in\sZ\times\sA} \int |c(x,a)-c(z,a)| \nu_{n,i_n(z)}(dx) \nonumber \\
&\leq 2 \omega_c(d_n). \label{compact:avebound2}
\end{align}
For the last term, we have
\begin{align}
|\hat{\rho}^n_{\hat{f}_n^{*}} - \rho_{f^{*}}| &= |\inf_{f\in\rF} \hat{\rho}^n_{f} - \inf_{f\in\rF} \rho_{f}| \leq \sup_{f\in\rF} |\hat{\rho}^n_{f} - \rho_{f}| \nonumber \\
&\leq 2 R \kappa^t \|c\| + \|c\| t \omega_p(d_n) \text{ (by (\ref{compact:avebound1}))}. \label{compact:avebound3}
\end{align}
Combining (\ref{compact:avebound1}), (\ref{compact:avebound2}), and (\ref{compact:avebound3}) implies
the result. \Halmos
\endproof

To explicitly calculate a convergence rate, we need to impose some structural assumptions on $\omega_c$ and $\omega_p$. One such assumption is linearity, which corresponds to the uniform Lipschitz continuity of $c(z,a)$ and $p(\,\cdot\,|z,a)$ in $z$. This means that $\omega_c(r) = K_1 r$ and $\omega_p(r) = K_2 r$, or equivalently, $|c(z,a)-c(y,a)| \leq K_1 d_{\sZ}(z,y)$ and $\|p(\,\cdot\,|z,a)-p(\,\cdot\,|y,a)\| \leq K_2 d_{\sZ}(z,y)$ for all $z,y \in \sZ$ and $a \in \sZ$. In this case, by Theorem~\ref{compact:mainthm4}, for all $t\geq1$ we have \begin{align}
|\rho_{\tilde{f}_n^{*}} - \rho_{f^{*}}| \leq 4 \|c\| R \kappa^t + 4 K_1  \alpha (1/n)^{1/d} + 4 \|c\| K_2 \alpha (1/n)^{1/d} t. \label{rateofconvergence}
\end{align}
To obtain a proper rate of convergence result (i.e., an upper bound that only depends on $n$) the dependence of the upper bound on $t$ has to be written as a function of $n$. This can be done by (approximately) minimizing the upper bound in (\ref{rateofconvergence}) with respect to $t$ for each $n$. Let us define the constants $I_1 \coloneqq 4 \|c\| R$, $I_2 \coloneqq 4 K_1 \alpha$, and $I_3 \coloneqq 4 \|c\| K_2 \alpha$. Then the upper bound in (\ref{rateofconvergence}) becomes
\begin{align}
I_1 \kappa^t + I_2 (1/n)^{1/d} + I_3 (1/n)^{1/d} t. \label{aux7}
\end{align}
For each $n$, it is straightforward to compute that
\begin{align}
t'(n) \coloneqq \ln\bigl(\frac{n^{1/d}}{I_4}\bigr) \frac{1}{ \ln(\frac{1}{\kappa})} \nonumber
\end{align}
is the zero of the derivative of the convex term in (\ref{aux7}), where $I_4 \coloneqq \frac{I_3}{I_1 \ln(\frac{1}{\kappa})}$. Letting $t = \lceil t'(n) \rceil$ in (\ref{aux7}), we obtain the following result.

\begin{corollary}\label{compact:corollary}
Suppose that $c(z,a)$ and $p(\,\cdot\,|z,a)$ are uniformly Lipschitz continuous in $z$ in addition to the assumptions imposed at the beginning of this section. Then, we have
\begin{align}
|\rho_{\tilde{f}_n^{*}} - \rho_{f^{*}}| \leq (I_1 I_4 + I_2) (1/n)^{1/d} + \frac{I_3}{\ln(1/\kappa)} (1/n)^{1/d} \ln\bigl(\frac{n^{1/d}}{I_4}\bigr). \nonumber
\end{align}
\end{corollary}

\section{Order Optimality for Approximation Errors in the Rate of Quantization.}\label{compact:order}

The following example demonstrates that the order of the performance losses in Theorem~\ref{compact:mainthm3} and Corollary \ref{compact:corollary} cannot be better than $O((\frac{1}{n})^{\frac{1}{d}})$. More precisely, we exhibit a simple standard example where we can lower bound the performance loss by $L(1/n)^{1/d}$, for some positive constant $L$. A similar result was obtained in \citet[Section IV]{SaLiYu13-2} for the case of quantization of action space, where the action space was a compact subset of $\R^m$ for some $m\geq1$. Therefore, when both state and action spaces are quantized, then the resulting construction is order optimal in the above sense as the approximation error, in this case, is bounded by the sum of the approximation errors in quantization of state space and quantization of action space.

In what follows $h(\,\cdot\,)$ and $h(\,\cdot\,|\,\cdot\,)$ denote differential and conditional differential entropies, respectively; see \citet[Chapter 8]{CoTh06}.

Consider the additive-noise system:
\begin{align}
z_{t+1} = F(z_t,a_t) + v_t, t=0,1,2,\ldots, \nonumber
\end{align}
where $z_t, a_t, v_t \in \R^d$. We assume that $\sup_{(z,a)\in\R^d\times\R^d} \frac{\|F(z,a)\|}{\|z\|+\|a\|}<1/2$.
The noise process $\{v_t\}$ is a sequence of i.i.d. random vectors whose common distribution has density $g$ supported on some compact subset $V$ of $\R^d$. We choose $V$ such that $\sZ = \sA$ can be taken to be compact subsets of $\R^d$. For simplicity suppose that the initial distribution $\mu$ has the same density $g$. It is assumed that the differential entropy $h(g) \coloneqq -\int_{\sZ} g(z) \log{g(z)} dz$ is finite. Let the one stage cost function be $c(z,a) \coloneqq \| z-a \|$. Clearly, the optimal stationary policy $f^{*}$ is induced by the identity $f^{*}(z) = z$, having the optimal cost $J(f^{*},\mu) = 0$ and $V(f^{*},\mu) = 0$. Let $\hf_n$ be the piece-wise constant extension of the optimal policy $f_n^{*}$ of the MDP$_n$ to the set $\sZ$. Fix $n\geq 1$ and define $D_t \coloneqq E_{\mu}^{\hf_n}\bigl[c(z_t,a_t)\bigr]$ for all $t$. Then, since $a_t = \hf_n(z_t)$ can take at most $n$ values in $\sA$, by the Shannon lower bound (SLB) (see \citet[p. 12]{YaTGr80}) we have for $t\geq1$
\begin{align}
\log{n} &\geq R(D_t) \geq h(z_t)+ \theta(D_t) \nonumber \\
&= h(F(z_{t-1},a_{t-1})+v_{t-1}) + \theta(D_t) \nonumber \\
&\geq h(F(z_{t-1},a_{t-1})+v_{t-1}|z_{t-1},a_{t-1}) + \theta(D_t) \label{compact:aux50} \\
&= h(v_{t-1}) + \theta(D_t), \label{compact:nneq25}
\end{align}
where $\theta(D_t) = - d + \log\biggl(\frac{1}{dV_d\Gamma(d)}\bigl(\frac{d}{D_t}\bigr)^d\biggr)$, $R(D_t)$ is the rate-distortion function of $z_t$, $V_d$ is the volume of the unit sphere $S_d = \{z: \|z\| \leq 1\}$, and $\Gamma$ is the gamma function. Here, (\ref{compact:aux50}) follows from the fact that conditioning reduces the entropy (see \citet[Theorem 2.6.5, p. 29]{CoTh06}) and (\ref{compact:nneq25}) follows from the independence of $v_{t-1}$ and the pair $(z_{t-1},a_{t-1})$. Note that $h(v_{t-1})=h(g)$ for all $t$. Thus, $D_t \geq L (1/n)^{1/d}$, where $L \coloneqq \frac{d}{2} \bigl(\frac{2^{h(g)}}{d V_d \Gamma(d)}\bigr)^{1/d}$. Since we have obtained stage-wise error bounds, these give
$|J(f^{*},\mu) - J(\hf_n,\mu)| \geq \frac{L}{1-\beta} (1/n)^{1/d} \text{  and  }
|V(f^{*},\mu) - V(\hf_n,\mu)| \geq L (1/n)^{1/d}$.

\vspace{3pt}

\begin{remark}
We note that if $h(x_{t+1}|x_t,a_t)$ can be lower bounded by some constant $k$ for all $t\geq1$, above analysis still holds by replacing $h(g)$ with $k$. For instance, this is the case if the transition probability $p(\,\cdot\,|x,a)$ admits a density which is bounded from above uniformly in $(x,a)$.
\end{remark}

\section{Numerical Examples.} \label{examples}

In this section, we consider two examples, the additive noise model and fisheries management problem, in order to illustrate our results numerically. Since computing true costs of the policies obtained from the finite models is intractable, we only compute the value functions of the finite models and illustrate their converge to the value function of the original MDP as $n\rightarrow\infty$.

Before proceeding to the examples, we note that all results in this paper apply  with straightforward modifications for the case of maximizing reward instead of minimizing cost.

\subsection{Additive Noise System.} \label{discounted example}

In this example, the additive noise system is given by
\begin{align}
x_{t+1}=F(x_{t},a_{t})+v_{t}, \text{ } t=0,1,2,\ldots \nonumber
\end{align}
where $x_t, a_t, v_t \in \R$ and $\sX = \R$. The noise process $\{v_{t}\}$ is a sequence of $\R$-valued i.i.d. random variables with common density $g$. Hence, the transition probability $p(\,\cdot\,|x,a)$ is given by
\begin{align}
p(D|x,a) = \int_{D} g(v-F(x,a)) m(dv) \text{  } \text{ for all $D\in \B(\R)$}, \nonumber
\end{align}
where $m$ is the Lebesgue measure. The one-stage cost function is $c(x,a) = (x-a)^2$, the action space is $\sA = [-L,L]$ for some $L>0$, and the cost function to be minimized is the discounted cost.

We assume that (i) $g$ is a Gaussian probability density function with zero mean and variance $\sigma^2$, (ii) $\sup_{a\in\sA} |F(x,a)|^2 \leq k_1 x^2 + k_2$ for some $k_1,k_2 \in \R_+$, (ii) $\beta < 1/\alpha$ for some $\alpha \geq k_1$, and (iv) $F$ is continuous. Hence, Assumption~\ref{as1} holds for this model with $w(x)=k+x^2$ and $M=4\bigl(\frac{L^2}{k}+x^2\bigr)$, for some $k \in \R_+$.

For the numerical results, we use the following parameters: $F(x,a)=x+a$, $\beta=0.3$,
$L=0.5$, and $\sigma = 0.1$.

We selected a sequence $\bigl\{[-l_n,l_n]\bigr\}_{n=1}^{15}$ of nested closed intervals, where $l_n = 0.5 + 0.25 n$, to approximate $\R$. Each interval is uniformly discretized using $\lceil 2 k_{\lceil\frac{n}{3}\rceil}l_n \rceil$ grid points, where $k_m = 5 m$ for $m=1,\ldots,5$ and $\lceil q \rceil$ denotes the smallest integer greater than or equal to $q \in \R$. Therefore, the discretization is gradually refined.
For each $n$, the finite state space is given by $\{x_{n,i}\}_{i=1}^{k_n} \cup \{\Delta_n\}$, where $\{x_{n,i}\}_{i=1}^{k_n}$ are the representation points in the uniform quantization of the closed interval $[-l_n,l_n]$ and $\Delta_n$ is a pseudo state. We also uniformly discretize the action space
$\sA = [-0.5,0.5]$ by using $2k_{\lceil\frac{n}{3}\rceil}$ grid points.
For each $n$, the finite state models are constructed as in Section~\ref{compact} by replacing $\sZ$ with $[-l_n,l_n]$ and by setting $\nu_n(\,\cdot\,) = \frac{1}{2} m_n(\,\cdot\,) + \frac{1}{2} \delta_{\Delta_n}(\,\cdot\,)$, where $m_n$ is the Lebesgue measure normalized over $[-l_n,l_n]$.

\begin{figure}[h]
\centering
\includegraphics[width=5in, height=2.5in]{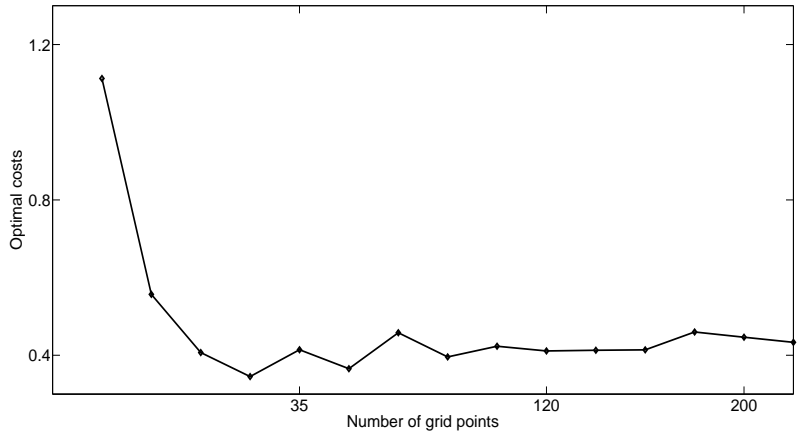}
\caption{Optimal costs of the finite models when the initial state is $x=0.7$}
\label{gr1}
\end{figure}

We use the value iteration algorithm to compute the value functions of the finite models. Figure~\ref{gr1} displays the graph of these value functions corresponding to the different values for the number of grid points, when the initial state is $x=0.7$.
The figure illustrates that the value functions of the finite models converge to the value function of the original model.

\subsection{Fisheries Management Problem.} \label{average example}

In this example we consider the following population growth model, called a Ricker model, see \citet[Section 1.3]{HeLa96}:
\begin{align}
x_{t+1} = \theta_1 a_t \exp\{-\theta_2 a_t + v_t\}, \text{ } t=0,1,2,\ldots \label{aux9}
\end{align}
where $\theta_1, \theta_2 \in \R_{+}$, $x_t$ is the population size in season $t$, and $a_t$ is the population to be left for spawning for the next season, or in other words, $x_t - a_t$ is the amount of fish captured in the season $t$. The one-stage `reward' function is $u(x_t-a_t)$, where $u$ is some utility function. In this model, the goal is to maximize the average reward.

The state and action spaces are $\sX=\sA=[\kappa_{\min},\kappa_{\max}]$, for some $\kappa_{\min}, \kappa_{\max} \in \R_{+}$. Since the population left for spawning cannot be greater than the total population, for each $x \in \sX$, the set of admissible actions is $\sA(x)=[\kappa_{\min},x]$ which is not consistent with our assumptions. However, we can (equivalently) reformulate above problem so that the admissible actions $\sA(x)$ will become $\sA$ for all $x\in\sX$. In this case, instead of dynamics in equation (\ref{aux9}) we have
\begin{align}
x_{t+1} = \theta_1 \min(a_t,x_t) \exp\{-\theta_2 \min(a_t,x_t) + v_t\}, \text{ } t=0,1,2,\ldots \nonumber
\end{align}
and $\sA(x) = [\kappa_{\min},\kappa_{\max}]$ for all $x\in\sX$. The one-stage reward function is $u(x_t-a_t)1_{\{x_t\geq a_t\}}$.

Since $\sX$ is already compact, it is sufficient to discretize $[\kappa_{\min},\kappa_{\max}]$. The noise process $\{v_{t}\}$ is a sequence of independent and identically distributed (i.i.d.) random variables which have common density $g$ supported on $[0,\lambda]$. Therefore, the transition probability $p(\,\cdot\,|x,a)$ is given by
\begin{align}
p\bigl(D|x,a\bigr) &= \Pr \biggl\{x_{t+1} \in D \biggl| x_t=x, a_t=a\biggr\} \nonumber \\
&= \Pr \biggl\{\theta_1 \min(a,x) \exp\{-\theta_2 \min(a,x) + v\} \in D\biggr\} \nonumber \\
&= \int_{D} g\biggl(\log(v) - \log(\theta_1 \min(a,x)) + \theta_2 \min(a,x)\biggr) \frac{1}{v} m(dv), \nonumber
\end{align}
for all $D\in \B(\R)$. To make the model consistent, we must have $\theta_1 y \exp\{-\theta_2 y + v\} \in [\kappa_{\min},\kappa_{\max}]$ for all $(y,v) \in [\kappa_{\min},\kappa_{\max}]\times[0,\lambda]$.

We assume that (i) $g > \epsilon$ for some $\epsilon \in \R_+$ on $[0,\lambda]$, (ii) $g$ is continuous on $[0,\lambda]$, and (iii) the utility function $u$ is continuous.
Define $h(v,x,a) \coloneqq g\bigl(\log(v) - \log(\theta_1 \min(a,x)) + \theta_2 \min(a,x)\bigr) \frac{1}{v}$, and for each $(x,a) \in \sX\times\sA$, let $S_{x,a}$ denote the support of $h(\,\cdot\,,x,a)$. Then, Assumption~\ref{compact:as2} holds for this model
with $\theta(x,a) = \inf_{v\in S_a} h(v,x,a)$ (provided that it is measurable), $\zeta = m_{\kappa}$ (Lebesgue measure restricted on $[\kappa_{\min},\kappa_{\max}]$), and for some $\lambda \in (0,1)$.

For the numerical results, we use the following values of the parameters:
\begin{align}
\theta_1 = 1.1, \text{ } \theta_2=0.1, \text{ }\kappa_{\max}=7, \text{ }\kappa_{\min}=0.005, \text{ }\lambda=0.5. \nonumber
\end{align}
We assume that the noise process is distributed uniformly over $[0,0.5]$. Hence, $g \equiv 1$ on $[0,0.5]$ and otherwise zero. The utility function $u$ is taken to be the shifted isoelastic utility function (see \citet[Section 4.1]{DuPr12})
\begin{align}
u(z) = 3 \bigl((z+0.5)^{1/3}-(0.5)^{1/3}\bigr). \nonumber
\end{align}
We selected 25 different values for the number $n$ of grid points to discretize the state space: $n=10,20,30,\ldots,250$.
The grid points are chosen uniformly over the interval $[\kappa_{\min},\kappa_{\max}]$. We also uniformly discretize the action space $\sA$ by using the following number of grid points: $5n= 50,100,150,\ldots,1250$.

We use the relative value iteration algorithm (see \citet[Chapter 4.3.1]{Ber95}) to compute the value functions of the finite models. For each $n$, the finite state models are constructed as in Section~\ref{compact} by replacing $\sZ$ with $[\kappa_{\min},\kappa_{\max}]$ and by setting $\nu_n(\,\cdot\,) = m_{\kappa}(\,\cdot\,)$.

\begin{figure}[h]
\centering
\includegraphics[width=5in, height=2.5in]{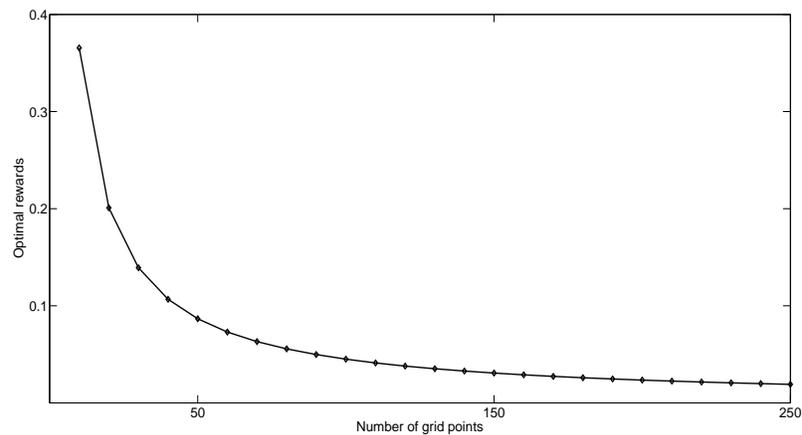}
\caption{Optimal rewards of the finite models when the initial state is $x=2$}
\label{gr3}
\end{figure}

Figure~\ref{gr3} shows the graph of the value functions of the finite models corresponding to the different values of $n$ (number of grid points), when the initial state is $x=2$. It can be seen that the value functions converge (to the value function of the original model).

\section{Conclusion.}\label{conc}

The approximation of a discrete time MDP by finite-state MDPs was considered for discounted and average costs for both compact and non-compact state spaces. Under usual conditions imposed for studying Markov decision processes, it was shown that if one uses a sufficiently large number of grid points to discretize the state space, then the resulting finite-state MDP yields a near optimal policy. Under the Lipschitz continuity of the transition probability and the one-stage cost function, explicit bounds were derived on the performance loss due to discretization in terms of the number of grid points for the compact state case. These results were then illustrated numerically by considering two different MDP models.

\section*{Acknowledgments.}
This research was supported in part by the Natural Sciences and Engineering Research Council (NSERC) of Canada. Parts of this work were presented at the American Control Conference in July 2015, Chicago, IL.

\bibliographystyle{ormsv080}
\bibliography{references}

\begin{thebibliography}{46}
\expandafter\ifx\csname natexlab\endcsname\relax\def\natexlab#1{#1}\fi
\expandafter\ifx\csname url\endcsname\relax
  \def\url#1{{\tt #1}}\fi
\expandafter\ifx\csname urlprefix\endcsname\relax\def\urlprefix{URL }\fi
\expandafter\ifx\csname urlstyle\endcsname\relax
  \expandafter\ifx\csname doi\endcsname\relax
  \def\doi#1{doi:\discretionary{}{}{}#1}\fi \else
  \expandafter\ifx\csname doi\endcsname\relax
  \def\doi{doi:\discretionary{}{}{}\begingroup \urlstyle{rm}\Url}\fi \fi

\bibitem[{Aliprantis and Border(2006)}]{AlBo06}
Aliprantis, C.D., K.C. Border. 2006.
\newblock {\it Infinite Dimensional Analysis\/}.
\newblock Springer.

\bibitem[{Bartoszynski(1961)}]{Bar61}
Bartoszynski, R. 1961.
\newblock A characterization of the weak convergence of measures.
\newblock {\it Ann. Math. Statist.\/} {\bf 32}(2) 561--576.

\bibitem[{Bertsekas and Shreve(1978)}]{BeSh78}
Bertsekas, D.~P., S.~E. Shreve. 1978.
\newblock {\it Stochastic optimal control: The discrete time case\/}.
\newblock Academic Press New York.

\bibitem[{Bertsekas(1975)}]{Ber75}
Bertsekas, D.P. 1975.
\newblock Convergence of discretization procedures in dynamic programming.
\newblock {\it IEEE Trans. Autom. Control\/} {\bf 20}(3) 415--419.

\bibitem[{Bertsekas(1995)}]{Ber95}
Bertsekas, D.P. 1995.
\newblock {\it Dynamic Programming and Optimal Control: Volume II\/}.
\newblock Athena Scientific.

\bibitem[{Bertsekas and Tsitsiklis(1996)}]{BeTs96}
Bertsekas, D.P., J.N. Tsitsiklis. 1996.
\newblock {\it Neuro-Dynammic Programming\/}.
\newblock Athena Scientific.

\bibitem[{Blackwell et~al.(1974)Blackwell, Freedman, and Orkin}]{BlFrOr74}
Blackwell, D., D.~Freedman, M.~Orkin. 1974.
\newblock The optimal reward operator in dynamic programming.
\newblock {\it Ann. Probab.\/} {\bf 2}(2) 926--941.

\bibitem[{Borkar(2002)}]{Bor02}
Borkar, V. 2002.
\newblock Convex analytic methods in {M}arkov decision processes.
\newblock E.A. Feinberg, A.~Shwartz, eds., {\it Handbook of Markov Decision
  Processes\/}. Kluwer Academic Publisher.

\bibitem[{Cavazos-Cadena(1986)}]{Cav86}
Cavazos-Cadena, R. 1986.
\newblock Finite-state approximations for denumerable state discounted {M}arkov
  decision processes.
\newblock {\it Appl. Math. Optim.\/} {\bf 14} 1--26.

\bibitem[{Chang et~al.(2007)Chang, Fu, Hu, and Marcus}]{ChFuHuMa07}
Chang, H.S., M.C. Fu, J.~Hu, S.I. Marcus. 2007.
\newblock A survey of some simulation-based methods in {M}arkov decision
  processes.
\newblock {\it Communications in Information System\/} {\bf 7} 59--92.

\bibitem[{Chow and Tsitsiklis(1991)}]{chow1991optimal}
Chow, C-S., J.~N. Tsitsiklis. 1991.
\newblock An optimal one-way multigrid algorithm for discrete-time stochastic
  control.
\newblock {\it IEEE Transactions on Automatic Control\/} {\bf 36}(8) 898--914.

\bibitem[{Cover and Thomas(2006)}]{CoTh06}
Cover, T.M., J.A. Thomas. 2006.
\newblock {\it Elements of Information Theory\/}.
\newblock 2nd ed. Wiley.

\bibitem[{Dufour and Prieto-Rumeau(2012)}]{DuPr12}
Dufour, F., T.~Prieto-Rumeau. 2012.
\newblock Approximation of {M}arkov decision processes with general state
  space.
\newblock {\it J. Math. Anal. Appl.\/} {\bf 388} 1254--1267.

\bibitem[{Dufour and Prieto-Rumeau(2013)}]{DuPr13}
Dufour, F., T.~Prieto-Rumeau. 2013.
\newblock Finite linear programming approximations of constrained discounted
  {M}arkov decision processes.
\newblock {\it SIAM J. Control Optim.\/} {\bf 51}(2) 1298--1324.

\bibitem[{Dufour and Prieto-Rumeau(2014)}]{DuPr14}
Dufour, F., T.~Prieto-Rumeau. 2014.
\newblock Approximation of average cost {M}arkov decision processes using
  empirical distributions and concentration inequalities.
\newblock {\it Stochastics\/}  1--35.

\bibitem[{Feinberg et~al.(2012)Feinberg, Kasyanov, and Zadioanchuk}]{FeKaZa12}
Feinberg, E.A., P.O. Kasyanov, N.V. Zadioanchuk. 2012.
\newblock Average cost {M}arkov decision processes with weakly continuous
  transition probabilities.
\newblock {\it Math. Oper. Res.\/} {\bf 37}(4) 591--607.

\bibitem[{Fox(1971)}]{Fox71}
Fox, B.L. 1971.
\newblock Finite-state approximations to denumerable state dynamic programs.
\newblock {\it J. Math. Anal. Appl.\/} {\bf 34} 665--670.

\bibitem[{Gordienko and Hernandez-Lerma(1995)}]{GoHe95}
Gordienko, E., O.~Hernandez-Lerma. 1995.
\newblock Average cost {M}arkov control processes with weighted norms:
  Existence of canonical policies.
\newblock {\it Appl. Math.\/} {\bf 23}(2) 199--218.

\bibitem[{Gray and Neuhoff(1998)}]{GrNe98}
Gray, G.M., D.L. Neuhoff. 1998.
\newblock Quantization.
\newblock {\it IEEE Trans. Inf. Theory\/} {\bf 44}(6) 2325--2383.

\bibitem[{Hern\'andez-Lerma(1989)}]{Her89}
Hern\'andez-Lerma, O. 1989.
\newblock {\it Adaptive {M}arkov Control Processes\/}.
\newblock Springer-Verlag.

\bibitem[{Hern\'andez-Lerma and Lasserre(1996)}]{HeLa96}
Hern\'andez-Lerma, O., J.B. Lasserre. 1996.
\newblock {\it Discrete-Time {M}arkov Control Processes: Basic Optimality
  Criteria\/}.
\newblock Springer.

\bibitem[{Hern\'andez-Lerma and Lasserre(1999)}]{HeLa99}
Hern\'andez-Lerma, O., J.B. Lasserre. 1999.
\newblock {\it Further Topics on Discrete-Time {M}arkov Control Processes\/}.
\newblock Springer.

\bibitem[{Hern\'andez-Lerma and Lasserre(2003)}]{HeLa03}
Hern\'andez-Lerma, O., J.B. Lasserre. 2003.
\newblock {\it {M}arkov Chains and Invariant Probabilities\/}.
\newblock Birkhauser.

\bibitem[{Hinderer(2005)}]{Hin05}
Hinderer, K. 2005.
\newblock Lipschitz continuity of value functions in {M}arkovian desision
  processes.
\newblock {\it Math. Meth. Oper. Res.\/} {\bf 62} 3--22.

\bibitem[{Jain and Varaiya(2006)}]{JaVa06}
Jain, R., P.P. Varaiya. 2006.
\newblock Simulation-based uniform value function estimates of {M}arkov
  decision processes.
\newblock {\it SIAM J. Control Optim.\/} {\bf 45}(5) 1633--1656.

\bibitem[{Ja\'{s}kiewicz and Nowak(2006)}]{JaNo06}
Ja\'{s}kiewicz, A., A.S. Nowak. 2006.
\newblock On the optimality equation for average cost {M}arkov control
  processes with {F}eller transition probabilities.
\newblock {\it J. Math. Anal. Appl.\/} {\bf 316} 495--509.

\bibitem[{Kuratowski(1966)}]{Kur66}
Kuratowski, K. 1966.
\newblock {\it Topology: Volume I\/}.
\newblock Academic Press Inc.

\bibitem[{Langen(1981)}]{Lan81}
Langen, H.J. 1981.
\newblock Convergence of dynamic programming models.
\newblock {\it Math. Oper. Res.\/} {\bf 6}(4) 493--512.

\bibitem[{Meyn and Tweedie(1993)}]{MeTw93}
Meyn, S.P., R.L. Tweedie. 1993.
\newblock {\it {M}arkov chains and stochastic stability\/}.
\newblock New York: Springer-Verlag.

\bibitem[{Ortner(2007)}]{Ort07}
Ortner, R. 2007.
\newblock Pseudometrics for state aggregation in average reward {M}arkov
  decision processes.
\newblock {\it Algorithmic Learning Theory\/}. Springer-Verlag.

\bibitem[{Puterman(2005)}]{Put05}
Puterman, M.L. 2005.
\newblock {\it {M}arkov Decision Processes\/}.
\newblock Wiley-Interscience.

\bibitem[{Ren and Krogh(2002)}]{ReKr02}
Ren, Z., B.H. Krogh. 2002.
\newblock State aggregation in {M}arkov decision processes.
\newblock {\it IEEE Conf. Decision Control\/}. Las Vegas, 3819 -- 3824.

\bibitem[{Saldi et~al.(2015)Saldi, Linder, and Y\"uksel}]{SaLiYu13-2}
Saldi, N., T.~Linder, S.~Y\"uksel. 2015.
\newblock Asymtotic optimality and rates of convergence of quantized stationary
  policies in stochastic control.
\newblock {\it IEEE Trans. Autom. Control\/} {\bf 60}(2) 553--558.

\bibitem[{Saldi et~al.(2016)Saldi, Y\"{u}ksel, and Linder}]{SaYuLi16}
Saldi, N., S.~Y\"{u}ksel, T.~Linder. 2016.
\newblock Near optimality of quantized policies in stochastic control under
  weak continuity conditions.
\newblock {\it J. Math. Anal. Appl.\/} {\bf 435} 321--337.

\bibitem[{Serfozo(1982)}]{Ser82}
Serfozo, R. 1982.
\newblock Convergence of {L}ebesgue integrals with varying measures.
\newblock {\it Sankhya Ser.A\/}  380--402.

\bibitem[{Shreve and Bertsekas(1979)}]{ShBe79}
Shreve, S.E., D.P. Bertsekas. 1979.
\newblock Universally measurable policies in dynamic programming.
\newblock {\it Math. Oper. Res.\/} {\bf 4}(1) 15--30.

\bibitem[{{Van Roy}(2006)}]{Roy06}
{Van Roy}, B. 2006.
\newblock Performance loss bounds for approximate value iteration with state
  aggregation.
\newblock {\it Math. Oper. Res.\/} {\bf 31}(2) 234--244.

\bibitem[{Vega-Amaya(2003)}]{Veg03}
Vega-Amaya, O. 2003.
\newblock The average cost optimality equation: a fixed point approach.
\newblock {\it Bol. Soc. Mat. Mexicana\/} {\bf 9}(3) 185--195.

\bibitem[{Villani(2009)}]{Vil09}
Villani, C. 2009.
\newblock {\it Optimal transport: old and new\/}.
\newblock Springer.

\bibitem[{White(1980)}]{Whi80}
White, D.J. 1980.
\newblock Finite-state approximations for denumerable state infinite horizon
  discounted {M}arkov decision processes.
\newblock {\it J. Math. Anal. Appl.\/} {\bf 74} 292--295.

\bibitem[{White(1982)}]{Whi82}
White, D.J. 1982.
\newblock Finite-state approximations for denumerable state infinite horizon
  discounted {M}arkov decision processes with unbounded rewards.
\newblock {\it J. Math. Anal. Appl.\/} {\bf 186} 292--306.

\bibitem[{Whitt(1978)}]{Whi78}
Whitt, W. 1978.
\newblock Approximations of dynamic programs {I}.
\newblock {\it Math. Oper. Res.\/} {\bf 3}(3) 231--243.

\bibitem[{Whitt(1979)}]{Whi79}
Whitt, W. 1979.
\newblock Approximations of dynamic programs {II}.
\newblock {\it Math. Oper. Res.\/} {\bf 4}(2) 179--185.

\bibitem[{Yamada et~al.(1980)Yamada, Tazaki, and Gray}]{YaTGr80}
Yamada, Y., S.~Tazaki, R.M. Gray. 1980.
\newblock Asymptotic performance of block quantizers with difference distortion
  measures.
\newblock {\it IEEE Trans. Inf. Theory\/} {\bf 26} 6--14.

\bibitem[{Yu and Bertsekas(2004)}]{YuBe04}
Yu, H., D.P. Bertsekas. 2004.
\newblock Discretized approximations for {POMDP} with average cost.
\newblock {\it The 20th Conf. UAI\/}. Banff, Canada.

\bibitem[{Y\"uksel and Ba\c{s}ar(2013)}]{YuBa13}
Y\"uksel, S., T.~Ba\c{s}ar. 2013.
\newblock {\it Stochastic Networked Control Systems: Stabilization and
  Optimization under Information Constraints\/}.
\newblock Boston, MA, Birkhauser.

\end{thebibliography}

\end{document}